\documentclass[a4paper,11pt]{article}
\usepackage{amsmath}
\usepackage{amssymb}
\usepackage[utf8]{inputenc}
\usepackage{dsfont}
\usepackage{authblk}
\usepackage{amsthm}
\usepackage{xcolor}
\usepackage{overpic}
\usepackage[
    inner=2.5cm,
    outer=2.5cm,
    top=2.7cm,
    ]
    {geometry}

\usepackage[colorlinks=true, allcolors=blue]{hyperref}

\numberwithin{equation}{section}

\def\Leb{{\mathrm{Leb}}}

\def\supp{{\mathrm{supp\;}}}

\providecommand{\keywords}[1]
{
  \small	
  \textbf{\textit{Keywords---}} #1
}

\theoremstyle{definition}
\newtheorem{defn}{Definition}[section]

\theoremstyle{plain}
\newtheorem{prop}[defn]{Proposition}
\newtheorem{lemma}[defn]{Lemma}

\newtheorem{coro}[defn]{Corollary}
\newtheorem{thmIntro}{Theorem}

\theoremstyle{remark}
\newtheorem{rmk}[defn]{Remark}

\usepackage[font=small]{caption}
\usepackage{subcaption}
\DeclareGraphicsExtensions{.png,.pdf,.jpg}

\def\supp{{\mathrm{supp}}}

\title{Tail behaviour of stationary densities\\ for one-dimensional random diffeomorphisms} 
\author{Jeroen S.W. Lamb, Guillermo Olic\'on-M\'endez\footnote{guillermo.olicon@uni-potsdam.de}, and Martin Rasmussen}

\begin{document}
\parskip 1.5mm
\maketitle

\begin{abstract}
	We study the asymptotic behaviour of stationary densities of one-dimensional random diffeomorphisms, at the boundaries of their support, which correspond to deterministic fixed points of extremal diffeomorphisms. In particular, we show how this stationary density at a boundary depends on the underlying noise distribution, as well as the linearisation of the extremal diffeomorphism at the boundary point (in case the corresponding fixed point is hyperbolic), or the leading nonlinear term of the extremal diffeomorphism (in case the corresponding fixed point is not hyperbolic).
\end{abstract}





\section{Introduction}

In the theory of Markov processes, the existence and uniqueness of stationary distributions is a classical and well-studied problem \cite{Doob53,Meyn93, Hairer11}. 
While stationary distributions depend intricately on details of the Markov process, their asymptotic behaviour at the boundary of their support may display universal asymptotic features. In this work, we consider discrete-time Markov processes generated by one-dimensional random diffeomorphisms, whose stationary densities are supported on a compact or semi-infinite interval, and discuss their asymptotic behaviour near the boundary points.

Tails of probability distributions, i.e.~the probability to be close to the boundary of their support,
describe the probability of a rare event occurring. In Markov processes, such rare events, represented by tails of stationary distributions, are often important in applications, as they may represent a specific risk. Stationary distributions of stochastic gradient flows with additive noise are Gibbs measures \cite{Jordan98} that admit explicit expressions in terms of the associated potential function, from which the tail behaviour as $x\rightarrow\pm\infty$
directly follows. In the discrete-time setting, the tails of the stationary distribution of the recurrence equation $x_{n+1}= a_n x_n + b_n$, where $(a_n,b_n)_{n\in\mathbb{Z}}$ are independent and identically distributed (i.i.d.) random variables, have  been studied \cite{Buraczewski16}. A nonlinear model, slightly more complicated than the aforementioned, was addressed in \cite{Borkovec01}, demonstrating power-law tails. In all the previous examples, the analysis is performed for distributions with unbounded support as $x\rightarrow\pm\infty$.

In this paper, we study tails of a stationary distribution $\mu$ and its density $\phi$ supported on an interval $M$ with $\sup M=x_+<\infty$, for a large class of Markov processes generated by one-dimensional random difference equations of the form $x_{n+1}=h(x_n,\omega_n)$, with 
$(\omega_n)_{n\in\mathbb{N}}$ representing i.i.d.~noise realisations 
that are sequences of random variables supported on $[-1,1]$. Such random dynamical systems with bounded noise naturally arise in many applications in different disciplines, see e.g.~\cite{Onofrio13} and references therein. Boundedness of the noise enables boundedness of the support of $\mu$, for example, when the system admits a compact \textit{minimal invariant set} \cite{Lamb15,Zmarrou07}. We show that if $h$ is sufficiently smooth and order-preserving for each fixed $\omega_0$, asymptotic
scaling laws of the tail exist as $x\rightarrow x_+$, that depend only on the noise distribution near its boundary point, that is $\omega_0=1$, and a key dynamical feature of the random map $h$. More specifically, $x_+$ is a fixed point of the so-called \textit{extremal map} $h_+(x):=\max_{\omega_0}h(x,\omega_0)$, whose Taylor expansion around $x_+$ is
\[
    h_+(x)=x_++\lambda(x-x_+)+\alpha(x_+-x)^r + o((x_+-x)^r),
\]
for some $\lambda\in(0,1]$ and integer $r\geq 2$. In case $\lambda<1$ (i.e.~$x_+$ is a hyperbolic fixed point of $h_+$), the tail behaviour of the stationary distribution depends on its Lyapunov exponent $\ln\lambda$. Alternatively, if $\lambda=1$ (i.e.~$x_+$ is a nonhyperbolic fixed point of $h_+$), it depends on $\alpha$ and $r$. 

Observe that the scaling law of $\phi$ near $x_+$ changes as $\lambda\rightarrow1$. This fact aligns well with the dynamical systems point of view,
where loss of hyperbolicity of a fixed point creates the possibility of a local bifurcation, see e.g.~\cite{Kuznetsov04}. Such a local bifurcation of $h_+$ is also known to have potential ramifications for the random difference equation as this can induce a discontinuous change in the support of the stationary distributions. This has been coined a \textit{topological bifurcation} in \cite{Lamb15}.
Therefore, the change in the scaling law anticipates this bifurcation. We illustrate the main results of this work in the following example.

\subsection{Example: an intermittency bifurcation}
	\label{SUBSEC: example}
\begin{figure}
     \centering
     \begin{subfigure}[b]{0.47\textwidth}
         \centering
          \begin{overpic}[width=\linewidth]{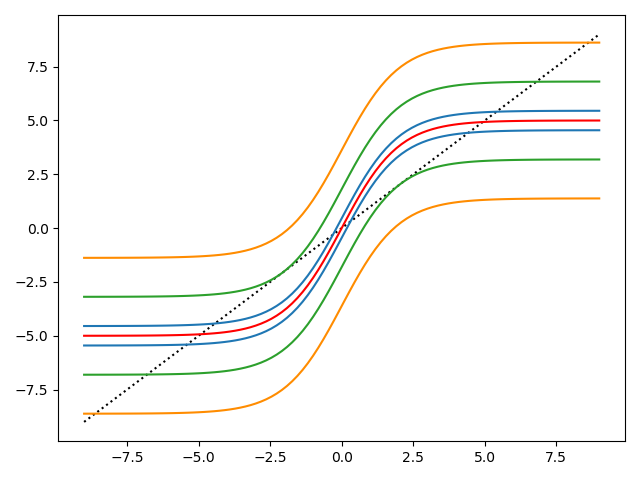}
           \put(50,-2){$x$}
      \put(12,55){$T(x)\pm \sigma$}
      \end{overpic}
      \caption{}
         \label{subfig:random map}
     \end{subfigure}
     \hfill
     \begin{subfigure}[b]{0.47\textwidth}
         \centering
          \begin{overpic}[width=\linewidth]{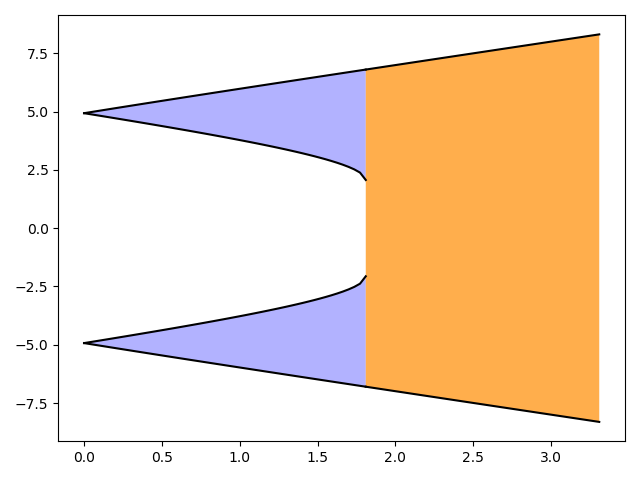}
           \put(50,-2){$\sigma$}
      \end{overpic}
         \caption{}
         \label{subfig:bifurcation}
     \end{subfigure}
     \hfill
       \caption[A topological bifurcation of minimal invariant sets]{In (a), a plot of the extremal maps $T(x)\pm\sigma$ of the system \eqref{eq:toy_example} for increasing values of $\sigma$ are shown, where the deterministic map is depicted in red (i.e. $\sigma=0$). Here $b=5$ was taken, and the blue, green, and orange plots correspond to noise strength values $\sigma=\sigma^*/4, \sigma^*$, and $2\sigma^*$, respectively, where $\sigma^*$ is given in \eqref{eq:sigma*}. As a reference, the identity line is portrayed as black dotted line.
       In (b), the change of the minimal invariant sets is presented as $\sigma$ increases, where the coexistence of two disjoint minimal invariant sets is portrayed in light blue colour, while in light orange the case when there is a unique one. The system exhibits a topological bifurcation at $\sigma=\sigma^*$.}
        \label{fig:bifMIS}
\end{figure}

 Consider the random map
	\begin{equation}
        \label{eq:toy_example}
		x_{n+1}=b\left(\frac{1-e^{-x}}{1+e^{-x}}\right)+\sigma\omega_n \equiv T(x_n)+\sigma\omega_n,
	\end{equation}
where $(\omega_n)_{n=0}^{\infty}$ is i.i.d.~sequence of uniformly distributed random variables on $[-1,1]$ and $b>2$ is a fixed parameter. We refer to the parameter $\sigma>0$ as the \textit{noise strength}. By changing the value of $\sigma$ we observe two different types of behaviours: there is a critical value $\sigma^*>0$ such that for $\sigma\leq \sigma^*$ there exist two stationary distributions supported on disjoint compact intervals. On the contrary, for $\sigma>\sigma^*$ there is a unique stationary distribution, see Figure~\ref{fig:bifMIS}. This type of topological bifurcation of invariant sets has been called an \textit{intermittency bifurcation} in \cite{Zmarrou07}. For the system \eqref{eq:toy_example}, the value $\sigma^*$ can be calculated explicitly by solving the system of equations
\[
    T'(x_+)=1, \qquad T(x_+)+\sigma=x_+,
\]
yielding in particular
\begin{equation}
    \label{eq:sigma*}
    \sigma^*=b\left( \frac{b-2+\sqrt{(b-1)^2-1}}{b+\sqrt{(b-1)^2-1}}\right) -\ln\left( b-1+\sqrt{(b-1)^2-1}\right).
\end{equation}

We illustrate the scaling laws of the upper tail of the stationary density (presented in more generality in Section~\ref{SUBSEC:setup_mainresults}) in the context of this example, comparing theoretical results with numerically obtained approximations, for which we employ \textit{Ulam's method} \cite{Ulam60}.

\begin{figure}
     \centering
     \begin{subfigure}[b]{0.47\textwidth}
         \centering
          \begin{overpic}[width=\linewidth]{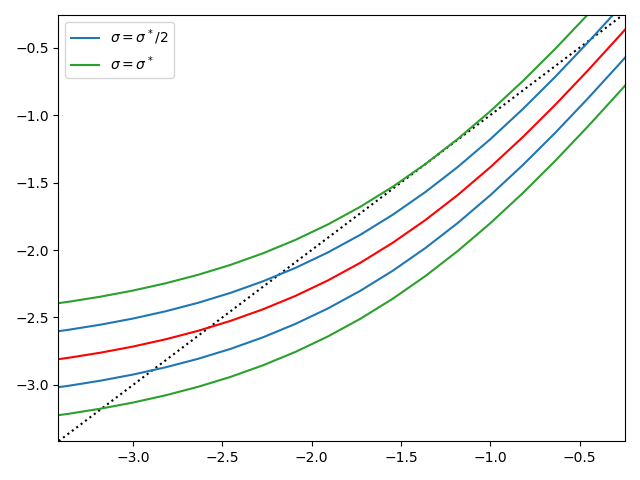}
            \put(50,-2){$x$}
      \put(12,55){$T(x)\pm \sigma$}
      \end{overpic}
      \caption{}
     \label{subfig:random map zoom}
     \end{subfigure}
     \hfill
     \begin{subfigure}[b]{0.47\textwidth}
         \centering
          \begin{overpic}[width=\linewidth]{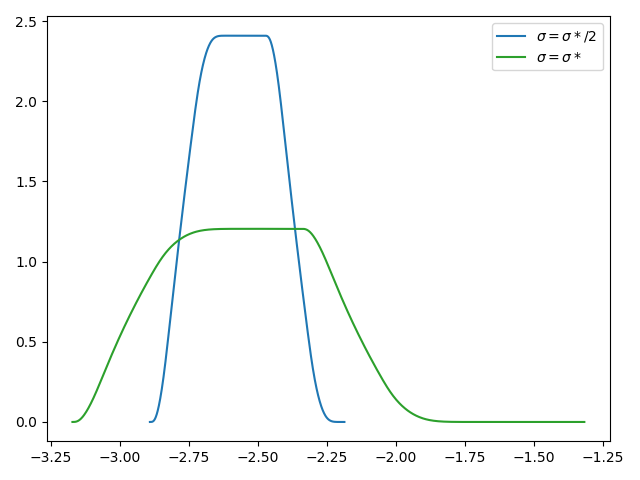}
            \put(50,-2){$x$}
      \put(-5,40){$\phi(x)$}
      \end{overpic}
      \caption{}
       
         \label{subfig:stationarydensities}
     \end{subfigure}
     \hfill
       \caption[Changes in stationary densities]{In (a), the minimal invariant set contained in $\mathbb{R}_-:=\{x\leq 0\}$ from Figure~\ref{fig:bifMIS} is depicted. Here $b=3$ was considered. The extremal maps are plotted for noise strengths $\sigma^*/2$ and $\sigma^*$, where $\sigma^*$ is the bifurcation value \eqref{eq:sigma*}. In (b), their corresponding stationary densities $\phi(x)$ are presented. The stationary density becomes flatter at the right boundary as $\sigma\rightarrow\sigma^*$.}
        \label{fig:stationary density}
\end{figure}

We consider the left invariant set for $\sigma<\sigma^*$. When $\sigma\rightarrow \sigma^*$ the shape of the stationary distribution becomes ``flatter'' at the boundary of its support, see
Figure~\ref{fig:stationary density}(b). 
 The boundary point $x_+$ is a hyperbolic fixed point of the extremal map $h_+(x):=T(x)+\sigma$.
 Theorem~\ref{THMasymptotichyperbolic} in Subsection~\ref{SUBSEC:setup_mainresults} then implies that 
the stationary density $\phi$ near $x_+$ satisfies 
\begin{equation}
\label{eq:example_hyp}
    \lim_{x\rightarrow x_+} \frac{\ln\phi(x)}{\ln^2(x_+-x)}=c_1:=\frac{1}{2\ln T'(x_+)}.,
\end{equation}
where, for simplicity, we denote $\ln^2 s \equiv \left(\ln s\right)^2$. When $\sigma=\sigma^*$, $x_+$ is a nonhyperbolic fixed point of $h_+$ (i.e. $T'(x_+)=1$) and the asymptotics of $\phi$, as given by Theorem~\ref{THMasymptoticNonhyperbolic} in Subsection~\ref{SUBSEC:setup_mainresults}, changes to
\begin{equation}
\label{eq:example_nonhyp}
    \lim_{x\rightarrow x_+} \frac{(x_+-x)\ln\phi(x)}{\ln(x_+-x)}=c_2:=\frac{4}{T''(x_+)}.
\end{equation}

\begin{figure}
     \centering
     \begin{subfigure}[b]{0.47\textwidth}
         \centering
          \begin{overpic}[width=\linewidth]{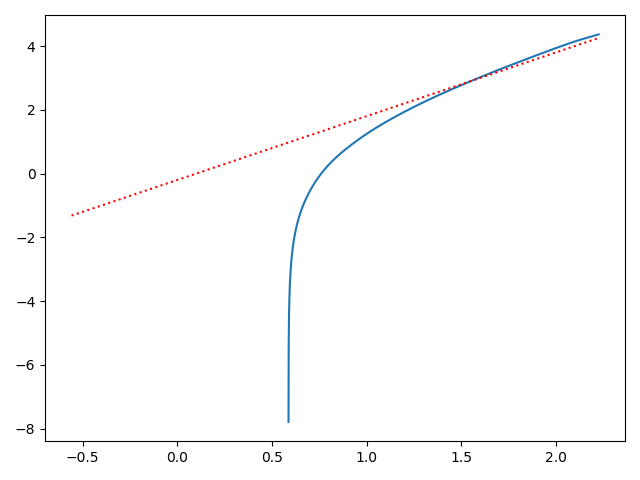}
           \put(60,-3){$\ln\ln\left( \frac{1}{x_+-x} \right)$}
      \put(60,45){$\ln\ln\left( \frac{1}{\phi(x)}\right)$}
      \put(15,65){$\ln\vert c_1\vert+2\ln\ln\left(\frac{1}{x_+-x} \right)$}
      \end{overpic}
      \caption{}
         \label{subfig:scale_hyperbolic}
     \end{subfigure}
     \hfill
     \begin{subfigure}[b]{0.47\textwidth}
         \centering
          \begin{overpic}[width=\linewidth]{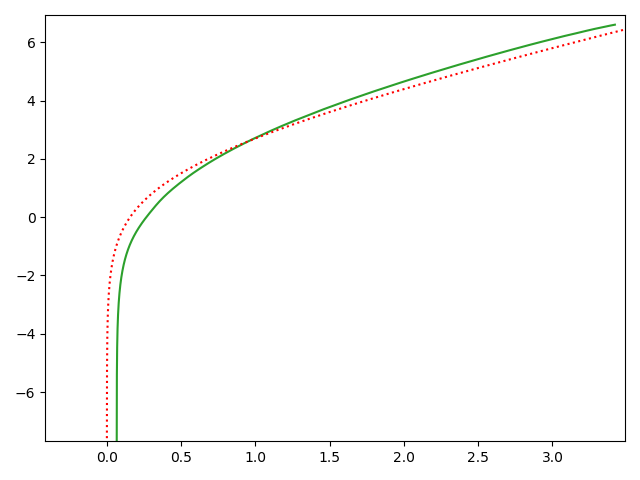}
          \put(60,-3){\small{$\ln\left( \frac{1}{x_+-x} \right)$}}
      \put(50,45){$\ln\ln\left( \frac{1}{\phi(x)}\right)$}
      \put(10,66){$\ln c_2+\ln\left(\frac{1}{x_+-x} \right)+\ln\ln\left(\frac{1}{x_+-x} \right)$}
      \end{overpic}
         \caption{}
         \label{subfig:scale_nonhyperbolic}
     \end{subfigure}
     \hfill
       \caption[Logarithm scaling for stationary densities]{
       A numerical validation of the main results for the system \eqref{eq:toy_example}. In (a), the value $\sigma=\sigma^*/2$ was taken, and thus the boundary $x_+$ is hyperbolic. We plot $\ln\ln (\phi(x))^{-1}$ as a function of $\ln \ln (x_+-x)^{-1}$ in blue, and compared with $\ln\vert c_1\vert+2\ln\ln\left(\frac{1}{x_+-x} \right)$ in red. In (b), the noise strength taken was $\sigma=\sigma^*$, so that $x_+$ is nonhyperbolic with $r=2$. The function $\ln\ln\left( \frac{1}{\phi(x)}\right)$ is plotted against $\ln\left( \frac{1}{x_+-x}\right)$ (in green), and compared to the function  $\ln c_2 +u + \ln u$ (in red) in the corresponding logarithmic scale.}
        \label{fig:scaling_densities}
\end{figure}

In order to facilitate the comparison between the theoretical results and numerical observations, we rewrite the scaling laws \eqref{eq:example_hyp} and \eqref{eq:example_nonhyp}. In the hyperbolic case one gets from \eqref{eq:example_hyp} that
\[
    \ln \ln\left( \frac{1}{\phi(x)} \right) = \ln(-c_1+o(1)) + 2\ln\ln\left(\frac{1}{x_+-x}\right)
\]
as $x\rightarrow x_+$, showing that the function $\ln(1/\phi(x))$, as a function of $\ln(1/(x_+-x))$, is a straight line in log-log scale with slope $2$ and $y$-intercept close to $\ln \vert c_1\vert$. In the nonhyperbolic case, we have from \eqref{eq:example_nonhyp} that, as $x\rightarrow x_+$,
\[
    \ln\ln\left(\frac{1}{\phi(x)}\right) = \ln(c_2+o(1)) + (r-1)\ln\left(\frac{1}{x_+-x}\right)+ \ln\ln\left(\frac{1}{x_+-x}\right).
\]
The function $\ln(1/\phi(x))$ as a function of $u=(x_+-x)^{-1}$ fits closely the function $u\mapsto \ln c_2+(r-1)u + \ln u$ in log-log scale. In Figure~\ref{fig:scaling_densities} we present a numerical validation in the hyperbolic and nonhyperbolic case\footnote{The reader can consult the Python script in \cite{OliconCode2024}.}.

\subsection{Set-up and main results}
\label{SUBSEC:setup_mainresults}
	
Throughout this work we consider the Markov process generated by the random difference equation
\begin{equation}
	\label{RDE}
		x_{n+1}=h(x_n,\omega_n),
\end{equation}
where $h:X\times\Delta\rightarrow X$ is a continuously differentiable function, defined on an interval $X\subset \mathbb{R}$ and $\Delta=[-1,1]$. The sequence $(\omega_n)_{n=0}^{\infty}$ is taken as an i.i.d.~sequence of random variables in $\Delta$, that is $(\omega_n)_{n=0}^{\infty}\in\Omega:=\Delta^{\mathbb{N}}$, where each $\omega_n$ is drawn according to a probability measure $\nu$ on the Borel space $(\Delta,\mathcal{B}(\Delta))$ with support $\supp \; \nu=\Delta$. We assume the following hypothesis.
\begin{itemize}
    \item[(H1)] \textbf{Monotonicity}.  $\partial_x h(x,\omega_0)>0$ and  $\partial_{\omega_0}h(x,\omega_0)>0$ for all $(x,\omega_0)\in X\times\Delta$.
    \item[(H2)] \textbf{Shape of the noise density}. The distribution $\nu$ on $(\Delta, \mathcal{B}(\Delta))$ is absolutely continuous with respect to the Lebesgue measure $\Leb$. Its density $p:\Delta\rightarrow \mathbb{R}_+=\{x\geq 0\}$ is a $C^k$ function for some $k\in\mathbb{N}_0$, such that $p(x)>0$ for every $x\in(-1,1)$ and for some $\beta>0$
    \begin{equation}
        \label{eq:noise_distribution}
        p(\omega_0)=\beta(1-\omega_0)^k+o(1-\omega_0)^k \quad \mbox{as } \omega_0\rightarrow1.
    \end{equation}
\end{itemize}
This setting includes uniformly distributed noise, where $p(\omega_0)\equiv 1/2$ and thus $k=0$.

There are two important implications of (H1). On the one hand, positivity of the partial derivative of $h$ with respect to $x$ implies that the map $h_{\omega_0}:=h(\cdot,\omega_0):X\rightarrow X$ is strictly monotonically increasing with a continuously differentiable inverse for each $\omega_0\in \Delta$. On the other hand, $\partial_{\omega_0}h>0$ implies that the set-valued map
\begin{equation}
    \label{eq:Set_valued_map}
        F(x):=h(x,\Delta)=
        [h_-(x),h_+(x)],
\end{equation}
where $h_{\pm}(x):=h(x,\pm1)$, is a nondegenerate compact interval for all $x\in X$. We refer to the maps $h_{+}$ and $h_-$ as the \textit{upper and lower extremal maps}, respectively.

As mentioned above, the extremal maps $h_{\pm}$ play a fundamental role in this paper, in part due to $x_+=\sup M$ being a fixed point of $h_+$, where an interval $M$ is the support of a stationary distribution. An analogous statement holds for $x_-:=\inf M>-\infty$. In particular, when the system remains confined to a compact set $X$, the supports of stationary distributions are compact minimal invariant intervals for the set-valued map $F$ \cite{Lamb15}. 

Stationary distributions of \eqref{RDE} are absolutely continuous with respect to the Lebesgue measure, and thus the stationary densities are fixed points of the transfer operator, see Section~\ref{SEC:preliminaries} below. The shape of these densities near $x_+$ strongly rely on iterating the explicit expression of the transfer operator $n$ times, where $n$  is intrinsically related to the speed at which $h_+^{n}(x)$ converges to $x_+$ as $n\rightarrow\infty$, and therefore inherently to the hyperbolicity of $x_+$, i.e. on whether $\lambda=h_+'(x_+)<1$ or not. Furthermore, in the case that $x_+$ is a nonhyperbolic fixed point of $h_+$, its shape depends on the first nonvanishing nonlinear term of $h_+$ in its Taylor expansion around $x_+$. In any case, as already pointed out in \cite{Zmarrou07}, the density of the stationary distribution is \textit{flat} at the boundaries  its support. However, as seen in the right boundary of the stationary densities in Figure~\ref{fig:stationary density} (b), in the hyperbolic case this density ``takes off'' quicker than in the nonhyperbolic scenario. Our main results provide a quantitative analysis for this phenomenon, as indicated in the following theorems.
\begin{thmIntro}[\textbf{Asymptotic scaling of $\phi$ near a hyperbolic boundary}]
\label{THMasymptotichyperbolic}
	Consider the system \eqref{RDE}, satisfying (H1) and (H2). Let $M$ be the support of a stationary distribution for \eqref{RDE}, with density $\phi:M\rightarrow\mathbb{R}_+$. If the upper extremal map $h_+$ is a $C^2$-(local) diffeomorphism, and $x_+:=\sup M<\infty$ is a hyperbolic fixed point of $h_+$ such that 
	\begin{equation}
		\label{hyperbolic}
			h_+(x)=x_++\lambda(x-x_+) +o(x_+-x) \quad \mbox{as } x\to x_+\,,
	\end{equation}
	for some $\lambda\in(0,1)$, then 
	\begin{equation}
		\label{orderhyperbolic}
			\lim_{x\rightarrow x_+}\frac{\ln\phi(x)}{\ln^2(x_+-x)} = \frac{k+1}{2\ln\lambda}\,.
	\end{equation}
\end{thmIntro}
It is clear from (\ref{orderhyperbolic}) that as $\lambda\rightarrow 1$, the limit (\ref{orderhyperbolic}) gives no further information, and a different scaling must be considered when $x_+$ is nonhyperbolic for $h_+$. We prove a corresponding result in such case.
\begin{thmIntro}[\textbf{Asymptotic scaling of $\phi$ near a nonhyperbolic boundary}]
\label{THMasymptoticNonhyperbolic}
	Consider the system \eqref{RDE}, satisfying (H1) and (H2). Let $M$ be the support of a stationary distribution for (\ref{RDE}), with density $\phi:M\rightarrow\mathbb{R}_+$. If the upper extremal map $h_+$ is a $C^{r+1}$-(local) diffeomorphism for some integer $r\geq 2$, and $x_+:=\sup M<\infty$ is a nonhyperbolic fixed point of $h_+$ such that 
	\begin{equation}
 \label{nonhyperbolic}
			h_+(x)=x+\alpha (x_+-x)^r+o((x_+-x)^r)\quad \mbox{as } x\to x_+\,,
	\end{equation}
	for some $\alpha>0$, then		\begin{equation}
		\label{ordernonhyperbolic}
			\lim_{x\rightarrow x_+}\frac{(x_+-x)^{r-1}\ln\phi(x)}{\ln(x_+-x)}= \frac{r(k+1)}{\alpha(r-1)}.
	\end{equation}
\end{thmIntro}

When $x_+$ is a nonhyperbolic fixed point with a quadratic contact to the identity line, i.e.~$r=2$, Theorem~\ref{THMasymptotichyperbolic} is not applicable and the scaling law to consider is given in \eqref{ordernonhyperbolic}. Analogously, if $\alpha\rightarrow 0$ for some $r\geq 2$ (that is when the system increases its degeneracy), (\ref{ordernonhyperbolic}) gives no further information and we need to analyse the same relationship but for the value $r+1$. 

In terms of the shape of the stationary density $\phi$, Theorem~\ref{THMasymptotichyperbolic} states that if $x_+$ is the right boundary of its support and this point is hyperbolic for the upper extremal map, then the stationary density $\phi$ can be approximated for $x\rightarrow x_+$  as
\begin{equation}
\label{eq:approx_hyp}
	 \phi(x) =(x_+-x)^{-(c_1+o(1))\ln(x_+-x)},
\end{equation} 
where $c_1\equiv c_1(k,\lambda):= \frac{k+1}{2\ln\lambda}<0$. On the other hand, when $x_+$ is a nonhyperbolic fixed point for the upper extremal map $h_+$ as indicated in Theorem~\ref{THMasymptoticNonhyperbolic}, as $x\rightarrow x_+$,
\begin{equation}
\label{eq:approx_nonhyp}
	\phi(x) = (x_+-x)^{\frac{c_2+o(1)}{(x_+-x)^{r-1}}},
\end{equation}
where $c_2\equiv c_2(k,r,\alpha):=\frac{r(k+1)}{\alpha(r-1)}>0$. These asymptotic expressions are \textit{universal} in the sense that they depend only on the local geometrical features of the extremal maps and the noise density.

\subsection{Outlook}
\label{SEC:conclusion}

Theorems \ref{THMasymptotichyperbolic} and \ref{THMasymptoticNonhyperbolic} establish a connection between the critical behaviour of the tails of stationary distributions and the dynamical features of the boundary points of their support, within the framework of random diffeomorphisms subjected to bounded noise. Specifically, the scaling laws depend on the Lyapunov exponent, in the case the boundary point is hyperbolic, or the first nonvanishing nonlinear term of the extremal map, for nonhyperbolic boundary points. 

It has previously been shown that the stationary densities are flat at the boundary of their support, even in higher dimensional systems \cite{Zmarrou07}, but no characterisation of their asymptotic behaviour was provided. To the best of our knowledge, our results are the first in this direction in the context of discrete-time Markov processes subjected to bounded noise. We remark, however, that some examples which relate bifurcation scenarios to the change of the shape of stationary densities were explored in one-dimensional SDEs, see \cite{Arnold92} and \cite[Section 9.3]{Arnold98}.

A straightforward extension of our results to order-reversing random diffeomorphisms is clear in certain circumstances. For instance, consider the i.i.d.~random difference equation
    $x_{n+1}=-\lambda x_{n}+\sigma\omega_n \equiv h(x_n,\omega_0)$,
where $\lambda\in(0,1)$ and the $\omega_i$ are distributed according to $d\nu = p\;dx$. In this case, we obtain an order preserving map by considering the random difference equation
 \[
    x_{n+1}=h^2(x_n,\tilde{\omega}_n) = \lambda^2 x_n + \sigma(\lambda\omega_{2n}+\omega_{2n+1})
 \]
 where $\tilde{\omega}_n:=\lambda\omega_{2n}+\omega_{2n+1}$. We expect that the transfer operator techniques developed in this work can be used for nonlinear order-reversing maps and more complicated frameworks.

In applications, it may be relevant to determine, based on real world data $(x_n)_{n=1}^{N}$, how close a system is to a bifurcation. For example, in some systems (like the one presented in Subsection~\ref{SUBSEC: example}), \textit{noise-induced tipping} \cite{Ashwin12} is exhibited only after a topological bifurcation has occurred. As already explored in this paper, in sufficiently smooth one-dimensional systems adhering to assumptions (H1) and (H2) this is indicated by $\lambda=h'_+(x_+)$, which is involved in the asymptotic scaling of the stationary density $\phi$ near $x_+$, cf. Theorem~\ref{THMasymptotichyperbolic}. However, two significant challenges arise. On the one hand, the value of $x_+$ is often unknown, and one can only rely on estimators like $Z_k=\max_{n\leq k} x_n$. On the other hand, the approximation of $\phi(x)$ typically requires a lot of data points near $x_+$, which are relatively inaccessible. In principle, one can explore whether the stationary distributions satisfy a \textit{law of rare events} \cite{Haan06} in order to estimate the characteristics of the tail.
A thorough exploration of this problem is beyond the scope of this paper and left for future work. 

The methods developed here rely on the features of the transfer operator near the boundary point. We expect them to be suitably adapted to different settings and variations from the problem addressed herein. For instance, a conditional version of the transfer operator can be used to study the density of the so-called \textit{quasi-stationary distributions}, and by this mean to study the asymptotic scaling of escape rates near a \textit{nondegenerate saddle-node bifurcation}, see \cite[Chapter 5]{Olicon21}. An improvement of these results follow in a similar way as those developed in this work, and will be communicated elsewhere. The analysis of quasi-stationary distributions becomes relevant in applications as a way to understand dynamical transient behaviour, which is of great importance, for instance, in ecological contexts \cite{Hastings21}. 

In the context of higher-dimensional systems, formulating scaling laws like Theorems \ref{THMasymptotichyperbolic} and \ref{THMasymptoticNonhyperbolic} remains as open problems. A first approach to this problem should be considered in two-dimensional systems, where the geometrical features of the boundary sets of $\varepsilon$-neighbourhoods of compact sets have already been studied \cite{Timperi20}. Moreover, the persistence of the minimal invariant sets in higher dimensions is characterised by the normal hyperbolicity of the so-called \textit{boundary map} \cite{TeyTimperi23}. We expect these novel techniques to be useful to translate topological bifurcations into changes in the asymptotic behaviour of stationary distributions.

\subsection*{Structure of the paper}
In Section~\ref{SEC:preliminaries} we provide the reader with some preliminaries used throughout this work, where stationary densities and the transfer operator \cite{Zmarrou07} are presented. In Section~\ref{SEC:main_lemma} we prove a technical lemma which provides the base upper and lower bounds used to prove the main theorems. The proofs of Theorems~\ref{THMasymptotichyperbolic} and \ref{THMasymptoticNonhyperbolic} are contained in Section~\ref{SECTIONtails}.

\section{Preliminaries}
\label{SEC:preliminaries}

We present a summary of the functional analytic tools employed throughout this paper. Theorems \ref{THMasymptotichyperbolic} and \ref{THMasymptoticNonhyperbolic} rely on the additional assumption that the system \eqref{RDE} admits a stationary distribution $\mu$. This is the case, for instance, when the random map admits a positively invariant compact set \cite{Araujo00,Kifer86, Zmarrou07}. The material we present here is based on \cite{Hurth22}.

\subsection{Notation}
We denote the Borel space $(M,\mathcal{B}(M))$ for any measurable set $M\subset \mathbb{R}$, and the space of probability measures on $(M,\mathcal{B}(M))$ as $\mathcal{M}_1(M)$. Furthermore, we denote the Lebesgue spaces $L^1(M)\equiv L^1(M,\mathcal{B}(M),\Leb)$.

Throughout this paper, when fixing $\omega_0\in \Delta$ or  $x\in X$, we denote the respective \textit{section maps} as
\[
    h_{\omega_0}:X\rightarrow X, \qquad h_x:\Delta\rightarrow X,
\]
defined simply as $h_x(\omega_0)=h(x,\omega_0)=h_{\omega_0}(x)$. Due to (H1), the section maps are $C^1$-diffeomorphisms on their images so that, for instance, $h\left(x,h_x^{-1}(y)\right)=y$. Moreover, since the inverse is continuously differentiable, by using the chain rule one gets
\begin{equation}
    \label{eq:derivative_inverse}
        \partial_x h_x^{-1}(y)=- \frac{\partial_x h\left(x,h_x^{-1}(y)\right)}{\partial_{\omega_0} h\left(x,h_x^{-1}(y)\right)}.
\end{equation}
In particular, we will need to compute $\omega_0=h_x^{-1}(y)$ for $x=h_+^{-1}(y)$. In other words, $\omega_0$ is such that $h\left(h_+^{-1}(y),\omega_0  \right)=y$,
which implies that
\begin{equation}
    \label{eq:inv_special}
    h_{h_+^{-1}(y)}^{-1}(y)=1.
\end{equation}

As indicated in Theorems \ref{THMasymptotichyperbolic} and \ref{THMasymptoticNonhyperbolic}, we study the asymptotic behaviour of different functions as $x\rightarrow x_+$, or equivalently as $x_+-x\rightarrow 0$. We say that $f(x)=O(g(x))$ as $x\rightarrow0$ if there is $C>0$ such that 
\[
    \vert f(x) \vert \leq C \vert g(x) \vert \quad \mbox{for all $x$ sufficiently close to $0$}\,.
\]
We say that $f(x)=o( g(x))$ as $x\rightarrow 0$ if
\[
    \lim_{x\rightarrow0} \frac{\vert f(x)\vert}{\vert g(x) \vert}=0.
\]
In the rest of this paper, these expressions will always refer to the asymptotic behaviour as the corresponding variable tending to $0$ and we omit this reference from here on.

\subsection{Stationary distributions and the transfer operator}
Recall that the system \eqref{RDE} induces a Markov chain in $X$. 
The \textit{transition probabilities} for the system (\ref{RDE}) are given as
\begin{equation}
	\label{transitionProbability}
 \mathbb{P}_x(A):=\int_X{\mathds{1}_A(h_{\omega_0}(x))d\nu(\omega_0)}.
\end{equation}
The family of probability distributions $\left( \mathbb{P}_x\right)_{x\in X}$ is also called the \emph{Markov kernel} of \eqref{RDE}. Notice that $\mathbb{P}_x$ is absolutely continuous with respect to $\Leb$ for each $x\in X$. Indeed, by doing the change of variable $y=h_{\omega_0}(x)$ we get that for all $A\in\mathcal{B}(X)$
\[
	\mathbb{P}_x(A)=\int_A{k(x,y)dy},
\]
where the transition density $k:X\times X\rightarrow \mathbb{R}_+$ is given by 
\begin{equation}
	\label{transitiondensity}
		k(x,y)=\frac{p\left( h_x^{-1}(y)\right)}{\partial_{\omega_0}h\left( x,h_x^{-1}(y) \right)}.
\end{equation}
Since the support of $p$ is $\Delta$, then $\supp(k(x,\cdot))=F(x)=\left[h_-(x),h_+(x)\right]$. In order to calculate the probability of landing on a measurable set $A$ at time $n$ when the initial condition is $x\in X$, we use the family of probability measures $\mathbb{P}^n_x$ defined as\footnote{Notice that, in fact, \eqref{transitionN} only depends on the first $n$ entries of $\omega$.}
\begin{equation}
	\label{transitionN}
 \mathbb{P}^n_x(A)= \int_X{\mathds{1}_A( h^n(x,\omega))d\mathbb{P}(\omega)},
\end{equation}
where $h^n(x,\omega):= h_{\omega_{n-1}}\circ\cdots\circ h_{\omega_0}(x)$. 
Again, the probability measures $\mathbb{P}_x^n$ are absolutely continuous with respect to $\Leb$ and $d\mathbb{P}^n_x=k_n(x,\cdot)\;d\Leb$, where $k_1\equiv k$ and
\begin{equation}
	\label{transitionDensityN}
		k_{n}(x,y)=\int_X{k_{n-1}(x,z)k(z,y)dz}=\int_X{k(x,z)k_{n-1}(z,y)dz}
\end{equation}
for $n\geq 2$. Notice that $\supp(k_n(x,\cdot))=F^n(x)\equiv [h_-^n(x),h_+^n(x)]$, where $h^n_{\pm}$ denotes the $n$-th iterate of the corresponding extremal map.
\begin{prop}
\label{PROPtransitiondensitiescontinuous}
	For each $n\geq 2$, the transition densities $k_n$ are continuous in $X\times X$.
\end{prop}
\begin{proof}
    By using the triangle inequality, this result follows from the fact that $k$ is uniformly bounded, and the maps $x\mapsto k(x,\cdot)$ and $y\mapsto k(\cdot,y)$ are
    continuous in $L^1(X)$.
\end{proof}

The evolution of $\vartheta\in\mathcal{M}_1(X)$ along the system \eqref{RDE} is induced by the operator $\mathcal{P}:\mathcal{M}_1(X)\rightarrow\mathcal{M}_1(X)$ defined as
\begin{equation}
    \label{Markov} \mathcal{P}\vartheta(A)=\int_X{\mathbb{P}_x(A)d\vartheta(x)}
\end{equation}
for all $A\in\mathcal{B}(X)$. Observe that $\mathcal{P}^n\vartheta=\int_X{\mathbb{P}_x^nd\vartheta(x)}$, and in particular $\mathbb{P}^n_x=\mathcal{P}^n\delta_x$, where $\delta_x$ is the Dirac delta distribution. Because of this, $\mathcal{P}^n\vartheta$ is interpreted as the evolution of $\vartheta$ under the system (\ref{RDE}) at time $n$. Given any $\vartheta\in\mathcal{M}_1(X)$, the distribution $\mathcal{P}\vartheta$ is absolutely continuous with respect to Lebesgue. Let us consider the space of densities
\[
	\mathcal{D}(X)=\lbrace g\in L^1(X)\;\vert\; g\geq 0, \ \Vert g\Vert_{L^1}=1 \rbrace.
\]
We define the \textit{transfer operator} $L:\mathcal{D}(X)\rightarrow \mathcal{D}(X)$ as the Radon-Nykodym derivative $Lg=d\mathcal{P}\vartheta/d\Leb$, where $d\vartheta=g\;d\Leb$. By linearity, the transfer operator can be extended to $L^1(X)$. Alternatively, $L:L^1(X)\rightarrow L^1(X)$ is the dual operator of the Markov operator $P:L^{\infty}(X)\rightarrow L^{\infty}(X)$ defined as
\[
	Pf(x)=\int_{\Delta}{f\circ h_{\omega_0}(x)d\nu(\omega_0)}.
\]
This operator is also known as the \textit{stochastic Koopman operator} \cite{Mezic19}.

The notion of a statistical equilibrium is given by the fixed points of $\mathcal{P}$. We say that $\mu\in\mathcal{M}_1(X)$ is a \textit{stationary distribution} if it is a fixed point of $\mathcal{P}$. Since $\mathcal{P}\mu$ is absolutely continuous with respect to $\Leb$,
we call its density $\phi\in\mathcal{D}(X)$ a \textit{stationary density}. Equivalently, $\phi$ is a stationary density if it is a fixed point of $L$. 

From now on, we assume that the system is restricted to an interval $M:=\textup{supp}\;\phi$, where $x_-=\inf M$ and $x_+:=\sup M<\infty$. Let us discuss in detail the action of the transfer operator $L:L^1(M)\rightarrow L^1(M)$, which is expressed as 
\begin{equation}
	\label{transferIntegralFOrmula}
		Lg(x)=\int_M{k(y,x)g(y)dy},
\end{equation}
see \cite{Zmarrou07}. Substituting (\ref{transitiondensity}) in the last expression, we obtain the explicit formula
\begin{equation}
	\label{transferExplicit}
		Lg(x)=\int_{h_+^{-1}(x)\wedge x_-}^{h_-^{-1}(x)\vee x_+}{\frac{p\left( h_y^{-1}(x)\right)}{\partial_{\omega_0} h\left( y,h_y^{-1}(x)\right)}g(y)dy}.
\end{equation}

Since we are interested in studying the behaviour of the stationary density $\phi:M\rightarrow\mathbb{R}_+$ near $x=x_+$, observe that whenever $x\in \left[ x_c,x_+ \right]$, where $x_c=h_+(x_-)\wedge h_-(x_+)$, we have
\begin{equation}
	\label{stationaryNearBoundary}
	\phi(x)=\int_{h_+^{-1}(x)}^{x_+}{\frac{p\left( h_y^{-1}(x)\right)}{\partial_{\omega_0}h(y,h_y^{-1}(x))}\phi(y)dy}.
\end{equation}
From the continuity and positivity of the integrand in \eqref{stationaryNearBoundary}, it follows that $x_+$ is a fixed point of $h_+$. Indeed, the integral on the right hand side attains a positive value if $h_+(x_+)<x_+$, which would contradict the fact that $\phi(x_+)=0$. Furthermore, $\phi$ is flat at $x_+$ \cite{Zmarrou07}.

\section{A fundamental lemma} 
\label{SEC:main_lemma}

From now on, we assume without loss of generality that $x_+=0$. If $x_+$ is hyperbolic for $h_+$, we have that
\begin{equation}
    \label{eq:Taylor_linear}
        h_+(x)=\lambda x +o(x),
\end{equation}
where $\lambda=h_+'(0)\equiv \partial_x h(0,1)$. In the case that $x_+$ is nonhyperbolic and $h_+$ is sufficiently smooth, there is $r\geq 2$ and $\alpha>0$ so that
\begin{equation}
    \label{eq:Taylor_Nonlinear}
        h_+(x)=x +\alpha \vert x\vert^r+o(\vert x\vert^r).
\end{equation}
Both cases will be analysed separately. We also consider the quantity $\gamma:=\partial_{\omega_0}h(0,1)>0$.

Let $\varepsilon>0$ be a sufficiently small constant. For each $\varepsilon$, let $x_0<0$ be sufficiently close to $0$, such that
	\begin{itemize}
		\item[C1:] for all $x\in[h_+^{-1}(x_0),0]$, formula \eqref{stationaryNearBoundary} holds,
        \item[C2:] for all $x\in[x_0,0]$ and $y\in[h_+^{-1}(x_0),0]$, the noise density $p$ satisfies
        \[
            (\beta-\varepsilon)\left(1-h_y^{-1}(x)\right)^k \leq p(h_y^{-1}(x)) \leq 
            (\beta+\varepsilon)\left(1-h_y^{-1}(x)\right)^k,
        \]
        \item[C3:] for all $x\in[x_0,0]$ and $y\in[h_+^{-1}(x_0),0]$, the random map satisfies
        \[
            \gamma-\varepsilon \leq 
            \partial_{\omega_0} h\left(y,h_y^{-1}(x)\right) \leq
            \gamma+\varepsilon, 
        \]
        \item[C4:] for all $x\in[x_0,0]$ and $y\in[h_+^{-1}(x_0),0]$, 
    \[
        \left( \frac{\lambda-\varepsilon}{\gamma} \right)\left( y-h_+^{-1}(x) \right) \leq 1-h_y^{-1}(x) \leq
        \left(\frac{\lambda+\varepsilon}{\gamma}\right) \left( y-h_+^{-1}(x) \right),
    \]

    \item[C5:] for all $x\in[h_{+}^{-1}(x_0),0]$,
    \[
    \lambda-\varepsilon \leq h_+'(x)\leq \lambda +\varepsilon.
    \]
	\end{itemize}  

Conditions C2--C5 are obtained from the local properties of the random map $h$ and the noise density $p$. In particular, C4 follows from Taylor's theorem applied to $h_y^{-1}(x)$ around $y=h_+^{-1}(x)$, cf. \eqref{eq:inv_special} and \eqref{eq:derivative_inverse}.

For each $x_0$ as above, we consider the following deterministic hitting time 
\begin{equation}
    \label{firsthitting}
        n_{x_0}^x=\min\{n\in\mathbb{N}_0 : h_+^n(x_0)\in (x,0]\}.
\end{equation}
This quantity becomes relevant in the proof of Theorems \ref{THMasymptotichyperbolic} and \ref{THMasymptoticNonhyperbolic}, since we iterate the transfer operator as many times possible so that conditions C1--C5 remain valid. As we see later, the scaling law of $n_{x_0}^x$, as $x\rightarrow 0$, depends on the hyperbolicity of $x_+$ (i.e $\lambda<1$), or the first nonvanishing nonlinear term in its Taylor expansion when $x_+$ is nonhyperbolic.

The following lemma constitutes the general groundwork for proving Theorem~\ref{THMasymptotichyperbolic} and Theorem~\ref{THMasymptoticNonhyperbolic}. For simplicity, we introduce the notation
\begin{equation}
    \label{eq:notation_aux}
        X_m=h_+^m(x), \quad Y_m=h_+^m(y), \quad Z_m=h_+^m(z) \qquad m\in\mathbb{Z}.
\end{equation}

\begin{lemma}
	\label{LEMMAcrucial}
		Consider the system (\ref{RDE}) under assumptions (H1) and (H2), restricted to the support $M$ of the stationary density $\phi$, with $\sup M=0$. Let $B=\sup_{x\in [h_+^{-1}(x_0),0 ]}\phi(x)$.
  Then, for each $\varepsilon>0$ there exist $x_0<0$ and constants $C_i\equiv C_i(\varepsilon)>0$ ($i=1,2$) bounded away from zero such that for each $x\in[x_0,0]$ and $n= n_{x_0}^x$, the following holds:
  \begin{itemize}
      \item[(a)] Hyperbolic case: under the assumptions of Theorem~\ref{THMasymptotichyperbolic},
		\begin{equation}
			\label{ineqMainLemma_Hyperbolic}
   \begin{split}
				\frac{(C_1 k!)^n}{[n(k+1)-1]!}
    \left( \frac{1}{\lambda+\varepsilon} \right)^{\frac{n(n-1)(k+1)}{2}}\int_{X_{-1}}^0 
    \left(y-X_{-1}\right)^{n(k+1)-1} \phi\left(  Y_{-n+1}\right)dy \leq \\
     \leq \phi(x) \leq
     \frac{(C_2 k!)^n}{[n(k+1)]!} 
      \left( \frac{1}{\lambda-\varepsilon} \right)^{\frac{n(n-1)(k+1)}{2}}
     \vert X_{-1}\vert^{n(k+1)} B. \end{split}
	\end{equation}
 \item[(b)] Nonhyperbolic case: under the assumption of Theorem~\ref{THMasymptoticNonhyperbolic}
 \begin{equation}
			\label{ineqMainLemma_Nonhyperbolic}
   \begin{split}
				\frac{(C_1 k!)^n}{[n(k+1)-1]!}\int_{X_{-1}}^0 
    \left(y-X_{-1}\right)^{n(k+1)-1} \phi\left(  Y_{-n+1}\right)dy \leq \\
     \leq \phi(x) \leq
     \frac{(C_2 k!)^n}{[n(k+1)]!} \prod_{m=1}^{n-1}\prod_{j=1}^m\left(  \frac{1}{h_+'\left(X_{-j-1} \right)} \right)^{k+1} \vert X_{-1}\vert^{n(k+1)}  B. \end{split}
	\end{equation}
  \end{itemize}
\end{lemma}

In the proof of Lemma~\ref{LEMMAcrucial}, the integrals
\[
    \int_{X_{-1}}^{Z_l} (y-X_{-1})^n (z-Y_{-l})^m
dy\]
need to be estimated. We summarise the necessary bounds in the following statement.

\begin{prop}
    \label{PROP_auxiliary}
    For every $m,n\in\mathbb{N}_0$, any $x\in [x_0,0]$, and $1\leq l\leq n_{x_0}^x$, whenever $z$ is such that $Z_l\geq X_{-1}$, the following inequalities hold:
    \begin{itemize}
        \item[(a)] Hyperbolic case: under the assumptions of Theorem~\ref{THMasymptotichyperbolic},
    \begin{equation}
        \label{eq:inequality_prop_hyp}
    \begin{split}
    \frac{m!\cdot n!}{(m+n+1)!} \left(\frac{1}{\lambda+\varepsilon}\right)^{lm} \left( Z_l-X_{-1} \right)^{n+m+1} \leq \\
   \leq \int_{X_{-1}}^{Z_l}\left(y-X_{-1}\right)^n \left( z-h_+^{-l}(y) \right)^m dy \leq \\
  \leq \frac{m!\cdot n!}{(m+n+1)!} \left(\frac{1}{\lambda-\varepsilon}\right)^{lm} \left( Z_l-X_{-1} \right)^{n+m+1}    \end{split}
    \end{equation}
    \item[(b)] Nonhyperbolic case: under the assumptions of Theorem~\ref{THMasymptoticNonhyperbolic},
     \begin{equation}
        \label{eq:inequality_prop_nonhyp}
    \begin{split}
    \frac{m!\cdot n!}{(m+n+1)!} \left(\prod_{j=1}^{l}\frac{1}{h_+'(Z_{j-1})}\right)^m \left( Z_l-X_{-1} \right)^{n+m+1} \leq \\ 
    \leq \int_{X_{-1}}^{Z_l}\left(y-X_{-1}\right)^n \left( z-h_+^{-l}(y) \right)^m dy \leq \\
    \leq \frac{m!\cdot n!}{(m+n+1)!} \left(\prod_{j=1}^{l}\frac{1}{h_+'(X_{-j-1})}\right)^m \left( Z_l-X_{-1} \right)^{n+m+1}    \end{split}
    \end{equation}
    \end{itemize}
\end{prop}
\begin{proof}
   We aim to estimate the integral
   \[
    I:=\int_{X_{-1}}^{Z_{l}}\left(y-X_{-1}\right)^n \left( z- Y_{-l} \right)^m dy.
   \]
For $m=0$, the result is straightforward. Assume then that $m\geq 1$. Using integration by parts, taking
\[
    \begin{split}
    u=&\left( z- Y_{-l} \right)^m, \\
    du=& -\frac{m(z-Y_{-l})^{m-1}}{h_+'(Y_{-1})\cdots h_+'(Y_{-l})}dy,
    \end{split}
    \qquad
    \begin{split}
    dv=&\left(y-X_{-1}\right)^n dy,\\
    v=& \frac{(y-X_{-1})^{n+1}}{n+1},
    \end{split}
\]
we obtain 
\[
    I = \frac{m}{n+1} \int_{X_{-1}}^{Z_l}
    \frac{\left(y-X_{-1}\right)^{n+1} \left( z- Y_{-l} \right)^{m-1} }{h_+'(Y_{-1})\cdots h_+'(Y_{-l})} dy.
\]
In the hyperbolic case, for any $y\geq X_{-1}$ we have that $Y_{-j}\geq X_{-1-j}\geq h_+^{-1}(x_0)$ for any $j\leq n_{x_0}^x$. Hence, using condition C5 so that $h_+'(Y_{-j})\geq \lambda-\varepsilon$, we obtain
\begin{equation}
    \label{eq:step_hyperbolic}
    I \leq \frac{m}{n+1}
    \left( \frac{1}{\lambda-\varepsilon} \right)^{l}
    \int_{X_{-1}}^{Z_l}
    \left(y-X_{-1}\right)^{n+1} \left( z-Y_{-l} \right)^{m-1} dy.
\end{equation}
On the other hand, in the nonhyperbolic case, $h_+'$ is increasing implying that $h_+'(Y_{-j})\geq h_+(X_{-1-j})$. Therefore, 
\begin{equation}
\label{eq:step_nonhyperbolic}
    I \leq \frac{m}{n+1}\left( \prod_{j=1}^l \frac{1}{h_+'(X_{-j-1})} \right)\int_{X_{-1}}^{Z_l}\left(y-X_{-1}\right)^{n+1} \left( z-Y_{-l} \right)^{m-1} dy.
\end{equation}
By iterating the procedure in each case, we obtain the upper bounds in \eqref{eq:inequality_prop_hyp} and \eqref{eq:inequality_prop_nonhyp}. The lower bounds follow in an analogous manner.
\end{proof}

Let us proceed with the proof the main lemma.

\begin{proof}[Proof of Lemma~\ref{LEMMAcrucial}]
We split the proof in two parts, each corresponding to the upper and lower bound respectively.

\noindent\textit{Part 1: upper bounds.} 

\noindent Let $x\in(x_0,0)$. Using \eqref{stationaryNearBoundary}, and conditions C2 and C3, we have that
\[
    \phi(x)\leq \frac{\beta+\varepsilon}{\gamma-\varepsilon}\int_{X_{-1}}^0\left(1-h_y^{-1}(x)\right)^{k}\phi(y)dy.
\]
Using condition C4, it follows that for any $x\in(x_0,0)$,
\begin{equation}
\label{eq:first_iterate_ineq}
    \phi(x)\leq C_2\int_{X_{-1}}^0(y-X_{-1})^k\phi(y)dy, \qquad C_2:=\frac{\beta+\varepsilon}{\gamma-\varepsilon} 
    \left( \frac{\lambda+\varepsilon}{\gamma}\right)^k = \frac{\beta\lambda^k}{\gamma^{k+1}}+O(\varepsilon).
\end{equation}

Assume that $x_0\leq X_{-1}$. Then, we can use \eqref{eq:first_iterate_ineq} again, obtaining
\[
    \phi(x)\leq C_2^2\int_{X_{-1}}^0\int_{Y_{-1}}^0(y-X_{-1})^k(z-Y_{-1})^k\phi(z)dzdy.
\]
Observe that the integration region can be rewritten as
\[
    \left\{(y,z) : X_{-1}\leq y\leq Z_1, \ X_{-2}\leq z \leq 0 \right\}.
\]
Therefore, by interchanging the order of integration, it yields
\begin{equation}
\label{eq:basic set}
    \phi(x)\leq C_2^2
    \int_{X_{-2}}^0\left(\int_{X_{-1}}^{Z_1}(y-X_{-1})^k(z-Y_{-1})^kdy \right) \phi(z)dz.
\end{equation}
Here we split the proof in the hyperbolic and nonhyperbolic case.

\medskip

\noindent\textit{(a) Hyperbolic case.} Using \eqref{eq:inequality_prop_hyp} from Proposition~\ref{PROP_auxiliary} for $n=m=k$ and $l=1$, it yields
\[
    \phi(x)\leq C_2^2\frac{(k!)^2}{(2k+1)!}\left( \frac{1}{\lambda-\varepsilon} \right)^k
    \int_{X_{-2}}^0(Z_{1}-X_{-1})^{2k+1} \phi(z)dz.
\]
By performing a change of variable $y=Z_1$, so that $Y_{-1}=z$ and $dy=h_+'(Y_{-1})dz$, it follows that
\[
    \phi(x)\leq C_2^2\frac{(k!)^2}{(2k+1)!}\left( \frac{1}{\lambda-\varepsilon} \right)^{k+1}
    \int_{X_{-1}}^0(y-X_{-1})^{2k+1} \phi(Y_{-1})dy.
\]

Let us assume now that $x_0\leq h_+^{-2}(x)$. In that case, we use \eqref{eq:first_iterate_ineq} for $\phi(Y_{-1})$, yielding
\[
        \phi(x)\leq C_2^3\frac{(k!)^2}{(2k+1)!}\left( \frac{1}{\lambda-\varepsilon} \right)^{k+1}
    \int_{X_{-1}}^0  \int_{Y_{-2}}^{0}
    (y-X_{-1})^{2k+1} (z-Y_{-2})^k \phi(z)dzdy.
\]
As before, the integration region is the set
\[
\left\{ (y,z) : X_{-1}\leq y\leq Z_2, \ X_{-3} \leq z\leq 0 \right\},
\]
so that by interchanging the order of integration one gets
\[
    \phi(x) \leq C_2^3\frac{(k!)^2}{(2k+1)!}\left( \frac{1}{\lambda-\varepsilon} \right)^{k+1}
    \int_{X_{-3}}^0  
    \left(\int_{X_{-1}}^{Z_2}
    (y-X_{-1})^{2k+1} (z-Y_{-2})^k dy \right) \phi(z)dz.
\]

By means of \eqref{eq:inequality_prop_hyp} in Proposition~\ref{PROP_auxiliary}, similarly as before, we obtain the bound
\[
    \phi(x)\leq \frac{(C_2k!)^3}{(3k+2)!} 
    \left( \frac{1}{\lambda-\varepsilon} \right)^{(1+2)(k+1)}\int_{X_{-1}}^0
    (y-X_{-1})^{3k+2}\phi(Y_{-2})dy.
\]
In general, if $x_0\leq h_+^{-n+1}(x)$, that is if $n\leq n_{x_0}^x$, one has
\begin{equation}
    \label{eq:general_iteration}
    \phi(x)\leq \frac{(C_2 k!)^n}{(n(k+1)-1)!}\left( \frac{1}{\lambda-\varepsilon} \right)^{\frac{n(n-1)(k+1)}{2}} 
    \int_{X_{-1}}^0(y-X_{-1})^{n(k+1)-1}\phi(Y_{-n+1}) dy.
\end{equation}
By taking $n=n_{x_0}^x$, and since $\phi$ is bounded on $[h_+^{-1}(x_0),0]$, the upper bound in \eqref{ineqMainLemma_Hyperbolic} holds.

\medskip

\noindent\textit{(b) Nonhyperbolic case.} We continue the general calculation obtained in \eqref{eq:basic set}, employing \eqref{eq:inequality_prop_nonhyp} instead. Similarly as before, and recalling that $h_+'$ is increasing, we obtain
\[
    \phi(x)\leq C_2^2\frac{(k!)^2}{(2k+1)!}\left( \frac{1}{h_+'(X_{-2})} \right)^{k+1}
    \int_{X_{-1}}^0(y-X_{-1})^{2k+1} \phi(Y_{-1})dy.
\]
By repeating the argument, for any $n\leq n_{x_0}^x$, we get the bound
\[
    \phi(x)\leq \frac{(C_2k!)^n}{(n(k+1)-1)!}\prod_{m=1}^{n-1}\prod_{j=1}^m\left(  \frac{1}{h_+'\left(X_{-j-1} \right)} \right)^{k+1} 
    \int_{X_{-1}}^0(y-X_{-1})^{n(k+1)-1} \phi(Y_{-n+1})dy.
\]
By considering $n=n_{x_0}^x$, and the fact that $\phi$ is bounded on $[h_+^{-1}(x_0),0]$, the upper bound in \eqref{ineqMainLemma_Nonhyperbolic} follows.

\medskip

\noindent\textit{Part 2: lower bounds.} 

\noindent We now prove that the lower bounds in \eqref{ineqMainLemma_Hyperbolic} and \eqref{ineqMainLemma_Nonhyperbolic} hold. Similarly as in the previous case, by means of conditions C2 and C3 we obtain that
\begin{equation}
\label{eq:lower_first_step}
    \phi(x)\geq C_1\int_{X_{-1}}^0(y-X_{-1})^k\phi(y)dy, \qquad C_1:=\frac{\beta-\varepsilon}{\gamma+\varepsilon} \left(\frac{\lambda-\varepsilon}{\gamma} \right)^k = \frac{\beta \lambda^k}{\gamma^{k+1}}-O(\varepsilon).
\end{equation}
Analogously, assuming that $x_0\leq X_{-1}$ we use \eqref{eq:lower_first_step} again, in order to get
\[
    \phi(x)\geq C_1^2\int_{X_{-1}}^0\int_{Y_{-1}}^0(y-X_{-1})^k(z-Y_{-1})^k\phi(z)dzdy.
\]
After interchanging the order of integration, we obtain
\begin{equation}
\label{eq:basic_set_lower}
    \phi(x)\geq C_1^2
    \int_{X_{-2}}^0\left(\int_{X_{-1}}^{Z_1}(y-X_{-1})^k(z-Y_{-1})^kdy \right) \phi(z)dz.
\end{equation}
The lower bounds follow similarly as the upper bounds. Notice, however, that in the nonhyperbolic case, the lower bound in \eqref{eq:inequality_prop_nonhyp} depends on the integration variable, and thus the double product cannot be taken out of the integral. After repeating $n$ times, we obtain
\[
    \phi(x)\geq\frac{(C_1 k!)^n}{[n(k+1)-1]!}\int_{X_{-1}}^0 \prod_{m=1}^{n-1}\prod_{j=1}^m 
    \left(  \frac{1}{h_+'(Y_{-j})}\right)^{k+1}
    \left(y-X_{-1}\right)^{n(k+1)-1} \phi\left(  Y_{-n+1}\right)dy.
\]
Since $h_+'(x)<1$ for $x$ sufficiently close to $0$, the double product is greater than one. The lower bound in \eqref{ineqMainLemma_Nonhyperbolic} follows.

\end{proof}

\section{Tails of stationary distributions}
\label{SECTIONtails}
In this section our goal is to prove Theorem~\ref{THMasymptotichyperbolic} and Theorem~\ref{THMasymptoticNonhyperbolic}, which we do in separate subsections. The general strategy involves deriving a scaling law of the first hitting times $n_{x_0}^x$ as $x\rightarrow 0$, and the use of Lemma~\ref{LEMMAcrucial}.

\subsection{Hyperbolic boundaries}
	In this subsection we assume that $x_+$ is a hyperbolic fixed point for the extremal map $h_+$. As before, we assume without loss of generality that $x_+=0$. Therefore,
 \[
		h_+(x)=\lambda x+o(x),
	\]
where $\lambda\in(0,1)$. The main strategy relies in the application of Lemma~\ref{LEMMAcrucial} for as many times as possible. This is indeed given by the hitting time $n_{x_0}^x$, as in \eqref{firsthitting}. We first analyse the behaviour of $n_{x_0}^x$ as $x\rightarrow 0$. 

\begin{lemma}[\textbf{Scaling law for $n_{x_0}^x$ near a hyperbolic point}] 
	\label{LEMMAscalingdetHyp}
		Let $h_+$ be a $C^2$-diffeomorphism as above. Then, for all $x_0\in M\setminus\lbrace 0\rbrace$,
		\begin{equation}
		\label{scalingHyp}
			\lim_{x\rightarrow 0}\frac{n_{x_0}^x}{\ln\vert x\vert}=\frac{1}{\ln\lambda}.
		\end{equation}
\end{lemma}
\begin{proof}
Since $0$ is a hyperbolic fixed point for $h_+$, and $h_+$ is a $C^2$-diffeomorphism, we can smoothly linearise it near $0$ with a near-identity map, see for instance \cite[Theorem 2.5]{Belitskii03}. More precisely, there is $\tilde{x}<0$ and a $C^2$-diffeomorphism $g(x)=x+\tilde{g}(x)$ such that $\tilde{g}(0)=\tilde{g}'(0)=0$, and
for all $x\in[\tilde{x},0]$ and $n\in\mathbb{N}_0$,
\[
     h_+^n(x)=g^{-1}\left(\lambda^n g(x)\right).
\]
Notice that $g$ is increasing and negative in $(\tilde{x},0)$, if $\tilde{x}$ is sufficiently close to $0$.

Let $x_0\in[\tilde{x},0)$, and denote $n\equiv n_{x_0}^x$ for $x>x_0$, i.e. $n$ is the unique integer such that
\begin{equation}
\label{eq:first_hitting_cond}
    h_+^{n-1}(x_0)<x, \qquad h_+^{n}(x_0)>x.
\end{equation}
From the first inequality, it follows using that $g$ is increasing that
\[
    \lambda^{n-1}g(x_0) < x + \tilde{g}(x) \quad \Rightarrow \quad \lambda^{n-1}\vert g(x_0)\vert > \vert x\vert - \tilde{g}(x).
\]
Taking logarithms on both sides yields
\[
    (n-1)\ln\lambda + \ln\vert g(x_0)\vert > \ln\vert x\vert + \ln \left( 1- \frac{\tilde{g}(x)}{\vert x\vert} \right).
\]
Dividing by $\ln\vert x\vert \ln \lambda >0$ on both sides we obtain
\[
    \frac{n-1}{\ln \vert x\vert} +o(1) > \frac{1}{\ln\lambda} +o(1),
\]
so that taking $\liminf_{x\rightarrow 0}$ implies 
\[
    \liminf_{x\rightarrow0} \frac{n}{\ln \vert x\vert} \geq \frac{1}{\ln\lambda}.
\]
The $\limsup$ follows analogously considering instead the second inequality in \eqref{eq:first_hitting_cond}.

For $x_0\leq \tilde{x}<x$, we have that
\[
	n_{\tilde{x}}^x\leq n_{x_0}^x \leq n_{x_0}^{\tilde{x}}+n_{\tilde{x}}^x,
\]
then \eqref{scalingHyp} holds for all $x_0\in M\setminus\lbrace 0\rbrace$.
\end{proof}	

Let $x_0\in M\setminus\{0\}$ and denote $n\equiv n_{x_0}^x$ as before. We show that $\ln(n(k+1))!=o(\ln^2\vert x\vert)$.
\begin{coro}
    \label{CORO:log_factorial_zero}
    \[
        \lim_{x\rightarrow 0} \frac{\ln(n(k+1))!}{\ln^2\vert x\vert}= 0
    \]
\end{coro}
\begin{proof}
We expand $\ln(n(k+1))!$ into a sum so that
\begin{align*}
0 \leq \ln(n(k+1))! & =\sum_{j=1}^{n(k+1)}\ln j \leq \int_{2}^{n(k+1)+1}\ln s\;ds \\
    & \leq[n(k+1)+1][\ln(n(k+1)+1)-1] -2(\ln2 -1) \\
    & = n(k+1)\left[ \ln n + O(1) \right] + \ln n+O(1).
\end{align*}
Notice that, dividing by $\ln^2\vert x\vert>0$ gives
\[
    \frac{\ln(n(k+1))!}{\ln^2\vert x\vert} \leq 
    \frac{n}{\ln \vert x\vert} (k+1) \left[ \frac{\ln n}{\ln \vert x\vert} +o(1) \right] + o(1).
\]
Observe that $\ln n /\ln\vert x\vert \rightarrow 0$ as $x\rightarrow 0$. This implies that $\limsup_{x\rightarrow 0} \frac{\ln(n(k+1))!}{\ln^2\vert x\vert}\leq 0$, and the result follows.
\end{proof}

We now proceed to prove Theorem~\ref{THMasymptotichyperbolic}. 
\begin{proof}[Proof of Theorem~\ref{THMasymptotichyperbolic}] We split the proof in two parts, for the upper and lower bounds in \eqref{orderhyperbolic}, respectively.

\noindent \textit{Part 1: upper bound in (\ref{orderhyperbolic}).} 

\noindent From the upper bound in \eqref{ineqMainLemma_Hyperbolic} we have that
\begin{align*}
 \ln \phi(x)\leq & \; n\ln\left( C_2\cdot k! \right) - \ln \;(n(k+1))! - \frac{n(n-1)(k+1)}{2}\ln(\lambda-\varepsilon) \\
  & \quad +n(k+1)\ln\vert X_{-1}\vert + \ln B.
\end{align*}
It follows from Lemma~\ref{LEMMAscalingdetHyp} that $n\ln (C_2\cdot k!)=o\left( \ln^2\vert x\vert \right)$. We also have that $\ln\vert X_{-1}\vert/\ln\vert x\vert=1+o(1)$. Therefore,
\begin{align*}
    \frac{\ln\phi(x)}{\ln^2\vert x\vert}\leq
    -\frac{\ln\;(n(k+1))!}{\ln^2\vert x\vert} - \frac{n(n-1)(k+1)\ln(\lambda-\varepsilon)}{2\ln^2\vert x\vert}+
    \frac{n(k+1)}{\ln\vert x\vert} + o(1).
\end{align*}
We use Corollary~\ref{CORO:log_factorial_zero} for the first term above, and Lemma~\ref{LEMMAscalingdetHyp} for the second and third terms, so that by taking $\limsup$ as $x\rightarrow 0$ we obtain
\begin{equation}
\label{eq:eq:almost_upperbound}
    \limsup_{x\rightarrow 0} \frac{\ln\phi(x)}{\ln^2\vert x\vert}\leq - 
    \frac{(k+1)\ln(\lambda-\varepsilon)}{2\ln^2\lambda}+
    \frac{(k+1)}{\ln\lambda}.
\end{equation}

Since \eqref{eq:eq:almost_upperbound} holds for all $\varepsilon$ sufficiently small. The upper bound in \eqref{orderhyperbolic} follows after taking $\varepsilon\rightarrow0$.

\medskip

\noindent \textit{Part 2: lower bound in (\ref{orderhyperbolic}).}  

\noindent We now use the lower bound in \eqref{ineqMainLemma_Hyperbolic}. First, we bound the integral expression from below by integrating on the subinterval $\left[X_{-1},x\right]$. Since $Y_{-n+1}\in [x_0, h_+(x_0)]$ when $y\in[X_{-1},x]$, we have that $\phi(Y_{-n+1})\geq D:=\inf_{x\in[x_0,h_+(x_0)]}\phi(x)>0$, so that
\begin{align*}
    \int_{X_{-1}}^0{(y-X_{-1})^{n(k+1)-1}\phi(Y_{-n+1})dy} & \geq \int_{X_{-1}}^x{(y-X_{-1})^{n(k+1)-1}\phi(Y_{-n+1})dy}\\ 
    & \geq D\int_{X_{-1}}^x(y-X_{-1})^{n(k+1)-1}dy \\
    &=\frac{D (x-X_{-1})^{n(k+1)}}{n(k+1)}.
\end{align*}

For each $\varepsilon>0$, assume without loss of generality that 
\begin{equation}
\label{estimate4}
	\zeta \vert x\vert < x-X_{-1}
\end{equation}
for all $x\in(x_0,0)$, where $\zeta\equiv \zeta_{\varepsilon}=\frac{1}{\lambda}-1-\varepsilon>0$. By means of \eqref{estimate4}, it follows that
\begin{equation}
\label{eq:lower_integral_hyper}
\int_{X_{-1}}^0{(y-X_{-1})^{n(k+1)-1}\phi(Y_{-n+1})dy} \geq \frac{ D\zeta^{n(k+1)}\vert x\vert^{n(k+1)}}{n(k+1)}.
\end{equation}

Bounding \eqref{ineqMainLemma_Hyperbolic} from below using \eqref{eq:lower_integral_hyper} yields
\[
    \phi(x) \geq \frac{D\tilde{C}_1^{n}}{[n(k+1)]!} \left( \frac{1}{\lambda +\varepsilon} \right)^{\frac{n(n-1)(k+1)}{2}} \vert x\vert^{n(k+1)} ,
\]
where $\tilde{C}_1:=k!\cdot C_1 \zeta^{k+1} $.
By taking natural logarithm, we have that
\begin{equation}
\begin{split}
    \ln\phi(x) & \geq n\ln\tilde{C}_1-
    \ln\; (n(k+1))!-
    \frac{n(n-1)(k+1)}{2}\ln(\lambda+\varepsilon) 
     + n(k+1)\ln \vert x\vert + \ln D
\end{split}
\end{equation}

By taking $\liminf$, after dividing by $\ln^2\vert x\vert$, it follows from Lemma~\ref{LEMMAscalingdetHyp} and Corollary~\ref{CORO:log_factorial_zero} that
\begin{equation}
\label{eq:almost_lowerbound}
    \liminf_{x\rightarrow 0} \frac{\ln\phi(x)}{\ln^2\vert x\vert}\geq 
    -\frac{(k+1) \ln(\lambda+\varepsilon)}{2\ln^2\lambda} + 
    \frac{k+1}{\ln\lambda},
\end{equation}
which holds for all $\varepsilon$ sufficiently small. Hence, the lower bound in \eqref{ineqMainLemma_Hyperbolic} holds after taking $\varepsilon\rightarrow 0$, and thus the proof of the theorem is complete.

\end{proof}

We formulate an analogous result to Theorem~\ref{THMasymptotichyperbolic} for the stationary distribution $d\mu=\phi d\Leb$. The tail of $\mu$, defined as $\mathcal{T}(x):=\mu([x,x_+])$, admits the same asymptotic behaviour as its density as we see in the next corollary

\begin{coro}
\label{CORO:tail_measure_hyperbolic}
	Assume that the hypothesis of Theorem~\ref{THMasymptotichyperbolic} hold. Then,
	\begin{equation}
 \label{eq:tail_asymptotic}
		\lim_{x\rightarrow x_+}\frac{\ln\mathcal{T}(x)}{\ln^2(x_+-x)}= \frac{k+1}{2\ln\lambda}.
	\end{equation}
\end{coro}
\begin{proof}
As before, we assume without loss of generality that $x_+=0$. Let $\varepsilon>0$ be arbitrarily small, and $x_0$ sufficiently close to $0$ such that for all $y\in[x_0,0]$,
\begin{equation}
\label{eq:exponential_bound}
    e^{(\gamma-\varepsilon)\ln^2\vert y \vert}\leq\phi(y)\leq e^{(\gamma+\varepsilon)\ln^2\vert y \vert},
\end{equation}
where $\gamma = \frac{k+1}{2\ln\lambda}<0$. Hence, for any $x>x_0$,
\[
    \mathcal{T}(x)=\int_{x}^0\phi(y)dy\leq 
    \int_x^0{e^{(\gamma+\varepsilon)\ln^2\vert y \vert}dy} \leq \vert x\vert e^{(\gamma+\varepsilon)\ln^2\vert x\vert},
\]
where the last inequality holds, since the integrand is a decreasing function. After taking $\ln$ on both sides we obtain that
\[
    \ln\mathcal{T}(x)\leq \ln\vert x\vert +  (\gamma+\varepsilon) \ln^2\vert x \vert.
\]
Dividing by $\ln ^2\vert x \vert$, and taking $\limsup_{x\rightarrow0}$, it yields
\[
    \limsup_{x\rightarrow 0}\frac{\ln\mathcal{T}(x)}{\ln^2\vert x\vert}\leq (\gamma+\varepsilon),
\]
for all $\varepsilon$ small enough. Hence, the upper bound in \eqref{eq:tail_asymptotic} holds.

For the lower bound, it follows from \eqref{eq:exponential_bound} that
\[
    \mathcal{T}(x)\geq \int_{x}^0{e^{(\gamma-\varepsilon)\ln^2\vert y \vert}dy} = \vert x\vert\left( \frac{1}{\vert x\vert}\int_x^{0} e^{(\gamma-\varepsilon) \ln^2\vert y \vert} \right),
\]
After taking logarithms, and using Jensen's inequality, we obtain
\[
    \ln\mathcal{T}(x)\geq \ln\vert x\vert +\frac{\gamma-\varepsilon}{\vert x\vert} \int_0^{\vert x\vert}\ln^2 s\; ds =\ln \vert x\vert + (\gamma-\varepsilon) \left( \ln^2\vert x\vert -2(\ln\vert x\vert-1)\right).
\]
The lower bound in \eqref{eq:tail_asymptotic} is obtained after dividing by $\ln^2\vert x \vert$, taking $\liminf_x\rightarrow0$, and letting $\varepsilon\rightarrow0$.
\end{proof}

\subsection{Nonhyperbolic boundaries}
We proceed now with the proof of Theorem~\ref{THMasymptoticNonhyperbolic}. We assume that the upper extremal map $h_+$ is a $C^{r+1}$-(local) diffeomorphism, where $x_+$ is a nonhyperbolic fixed point (assumed to be $0$, without loss of generality), such that for some $\alpha>0$,
\begin{equation}
\label{MapNonHyperbolic}
	h_+(x)=x+\alpha \vert x\vert^r+o(x^r).
\end{equation}
Whenever $r=2$, we say that $0$ is a nondegenerate nonhyperbolic fixed point. 

Analogously to Theorem~\ref{THMasymptotichyperbolic}, we derive an asymptotic scaling of the first hitting time
\[
    n_{x_0}^x=\min\{ n\in\mathbb{N} : h_+^{n}(x_0)\in(x,0) \}.
\]
\begin{lemma}[\textbf{Scaling law for $n_{x_0}^x$ near a nonhyperbolic point}]
\label{LEMMAscalingdeterministic}
	Let $h_+$ be a $C^{r+1}$-diffeomorphism as above. For all $x_0\in M\setminus\lbrace 0\rbrace$,
	\begin{equation}
		\label{scalingLawNonHyperbolic}
			\lim_{x\rightarrow 0}\vert x\vert^{r-1}n_{x_0}^x=\frac{1}{\alpha(r-1)}
	\end{equation}
\end{lemma}
\begin{proof}
	We split the proof in three steps. First, we argue how to (locally) embed the map $h_+$ into a semiflow $\left(\psi_t\right)_{t\geq 0}$ near $0$, i.e. we find a continuous-time dynamical system whose time-one map $\psi_1$ coincides with $h_+$ for all $x\in(x_1,0]$, for some $x_1<0$. By doing so we give an expression of the \textit{ordinary differential equation} (ODE), for which the flow $\psi_t$ is a solution. Secondly, we derive a scaling law for the first-hitting time of the flow $\psi_t$, that is
	\[
		t_{x_0}^x:=\inf{\lbrace t\geq0 \;\vert\; \psi_t(x_0)>x \rbrace},
	\]
	as $x\rightarrow 0$. Last, we conclude that this scaling law coincides with that one for $n_{x_0}^x$ as given in (\ref{scalingLawNonHyperbolic}) for all $x_0\in M\setminus\{0\}$. For finer details the reader can consult Proposition 3.3, Lemma 3.4, and Corollary 3.5 in \cite{Olicon21}.

\medskip

\noindent \textit{Step 1: embedding $h_+$ into a flow $\psi_t$.} 

\noindent Due to \textit{Takens's embedding theorem}, there exist $x_1<0$ and a $C^{r+1}$ flow $\left(\psi_t\right)_{t\geq 0}$ such that $\psi_1(x)=h_+(x)$ for all $x\in(x_1,0]$, see \cite[Appendix 3]{Yoccoz95}. The flow $\psi_t$ satisfies an ODE
\[
	\dot{z}=G(z),
\]
where the vector field $G:(x_1,0]\rightarrow \mathbb{R}$ is $C^r$. Therefore, assume that
\[
	G(x)=a_1x+a_2x^2+\cdots+a_rx^r+o(x^r),
\]
\[
	\psi_t(x)=b_1(t)x+\cdots b_r(t) x^r+o(x^r).
\]
Since $\dot{\psi_t}(x)=G(\psi_t(x))$, by substituting both expressions and considering that
\[
\begin{split}
    b_1(0)=b_1(1)=1, \qquad  b_r(0)=0; \; b_r(1)=(-1)^r\alpha, \\
    b_j(0)=b_j(1)=0 \quad (j=2,\ldots,r \textup{ if } r>2).
\end{split}
\]
we obtain
\begin{equation}
	\label{vectorfieldexpression}
		G(x)=\alpha \vert x\vert^r+o(x^r).
\end{equation}

\noindent \textit{Step 2: first-hitting times for $\psi_t$ near $0$.}

\noindent Since $\psi_t(x)$ satisfies $\dot{z}=G(z)$, then given $x_0<x<0$ we have that
\begin{equation}
\label{eq:continuoustime_hitting}
	t_{x_0}^x=\int_{x_0}^x{\frac{dz}{G(z)}}.
\end{equation}
Let $\varepsilon\in(0,\alpha)$ be arbitrarily small, and consider $z_0\equiv z_0(\varepsilon) \in(x_1,0)$ such that $z_0\rightarrow0$ as $\varepsilon\rightarrow0$, and
\[
	(\alpha-\varepsilon)\vert z\vert^r\leq G(z)\leq (\alpha+\varepsilon)\vert z\vert^r
\]
for all $z\in[z_0,0]$. Using these estimates in \eqref{eq:continuoustime_hitting}, we first estimate $t_{z_0}^x$ as
\[
	\frac{(-1)^{r+1}}{(\alpha+\varepsilon)(r-1)}\left[\frac{z_0^{r-1}-x^{r-1}}{z_0^{r-1}x^{r-1}} \right]    \leq t_{z_0}^x\leq \frac{(-1)^{r+1}}{(\alpha-\varepsilon)(r-1)}\left[\frac{z_0^{r-1}-x^{r-1}}{z_0^{r-1}x^{r-1}} \right].
\]
Therefore, we have
\begin{equation}
\label{eq:limsup_limninf_hitting}
	\frac{1}{(\alpha+\varepsilon)(r-1)}\leq \liminf_{x\rightarrow 0}\vert x\vert^{r-1}t_{z_0}^x\leq \limsup_{x\rightarrow 0}\vert x\vert^{r-1}t_{z_0}^x\leq \frac{1}{(\alpha-\varepsilon)(r-1)}.
\end{equation}

For any $x_0\in(x_1,0)$, let $\varepsilon$ be arbitrarily small such that $x_0<z_0$. Since $t\mapsto\psi_t(x_0)$ is increasing, it is straightforward that $t_{x_0}^x=t_{x_0}^{z_0}+t_{z_0}^x$. Hence, after multiplying times $\vert x\vert^{r-1}$ and taking $\limsup$ and $\liminf$ we have that \eqref{eq:limsup_limninf_hitting} holds by replacing $z_0$ with $x_0$. Since $x_0$ is independent of $\varepsilon$, we obtain that
\begin{equation}
	\label{scalinglawContinuousTime}
	\lim_{x\rightarrow 0}\vert x\vert^{r-1}\cdot t_{x_0}^x=\frac{1}{\alpha(r-1)}.
\end{equation}

\noindent \textit{Step 3: first-hitting time for $h_+$ near $x=0$.} 

\noindent We show that (\ref{scalingLawNonHyperbolic}) holds. First, let us consider $x_0\in(x_1,0)$ for $x_1$ as given in Step 2, so that
\[
	t_{x_0}^x \leq n_{x_0}^x\leq t_{x_0}^x+1,
\]
since $h_+$ is embedded in $\psi_t$ in $(x_1,0]$. From Step 2 above, the result follows for $x\in(x_1,0)$ by multiplying times $\vert x\vert^{r-1}$ and taking limit as $x\rightarrow 0$. The extension to any $x_0\in M\setminus\{0\}$ follows the same lines as in Lemma~\ref{LEMMAscalingdetHyp}.
\end{proof}

\begin{rmk}
	For an arbitrary $x_0\in (x_1,0)$ the map $x\longmapsto t_{x_0}^x$ is a diffeomorphism from $(x_0,0)$ to $(0,\infty)$ which is strictly increasing. Its inverse is precisely $t\longmapsto \psi_t(x_0)$. Therefore, by means of a change of variables, it is true that
	\[
		\lim_{t\rightarrow\infty}\vert\psi_t(x)\vert^{r-1}\cdot t=\frac{1}{\alpha(r-1)}
	\]
for all $x\in (x_1,0)$. Since $h_+$ is embedded in $\psi_t$, by taking the sequence $t_n=n$ for $n\in\mathbb{N}$ we obtain that
\begin{equation}
\label{scalingInverse}
	\lim_{n\rightarrow\infty}\vert h_+^n(x)\vert^{r-1}\cdot n=\frac{1}{\alpha(r-1)}.
\end{equation}
In a similar way as before, we can extend the limit (\ref{scalingInverse}) for all $x\in M\setminus\lbrace0\rbrace$. This way we have obtained an equivalent formulation of Lemma 2.2. in \cite{Huls01}. This result shows the close relationship between the speed of convergence towards nonhyperbolic fixed points and the hitting times of shrinking neighbourhoods around such fixed points.
\end{rmk}

Similarly as in the hyperbolic case, we need to deal with the term $\ln n(k+1)!$, where $n\equiv n_{x_0}^x$, for which we have the following result.
\begin{coro}
    \label{CORO: log_factorial_nonhyp}
    \[
        \lim_{x\rightarrow0} \frac{\vert x\vert^{r-1}\ln n(k+1)!}{\ln\vert x\vert} = -\frac{k+1}{\alpha}.
    \]
\end{coro}
\begin{proof}
First, we show that for all $x_0\neq 0$,
\begin{equation}
    \label{eq:first_limit}
   \lim_{x\rightarrow 0}\frac{\ln n}{\ln \vert x\vert}=-(r-1).
\end{equation}
Indeed, let $x_0\neq 0$ and $\varepsilon>0$ sufficiently small. From Lemma~\ref{LEMMAscalingdeterministic}, for every $x$ close enough to $0$,
\[
    \frac{1}{\alpha(r-1)}-\varepsilon\leq \vert x\vert^{r-1}\cdot n_{x_0}^x\leq \frac{1}{\alpha(r-1)}+\varepsilon.
\]
By taking $\ln$, and dividing by $\ln\vert x\vert <0$, it yields
\[
     \frac{\ln\left(  \frac{1}{\alpha(r-1)}+\varepsilon \right)}{\ln\vert x\vert}
    \leq(r-1)+\frac{\ln n_{x_0}^x}{\ln\vert x\vert}\leq\
    \frac{\ln\left(  \frac{1}{\alpha(r-1)}-\varepsilon \right)}{\ln\vert x\vert}.
\]
The claim follows by taking limits as $x\rightarrow0$.

    Now, similarly as in Corollary~\ref{CORO:log_factorial_zero},
    \begin{align*}
    \ln(n(k+1))! & \leq  \int_2^{n(k+1)+1}\ln s\;ds = [n(k+1)+1][\ln (n(k+1)+1)-1]-2(\ln2-1) \\
    & =  n(k+1) \left[ \ln n +O(1) \right] + \ln n + O(1)
    \end{align*}
Multiplying both sides bx $\vert x\vert^{r-1}/\ln\vert x\vert <0$ yields
\[
    \frac{\vert x\vert^{r-1}\ln (n(k+1))!}{\ln\vert x\vert} \geq  (k+1) \cdot \left( \vert x\vert^{r-1} n \right) \cdot \left( \frac{\ln n}{\ln \vert x\vert} + o(1)\right ) + o(1)
\]

Due to Lemma~\ref{LEMMAscalingdeterministic} and \eqref{eq:first_limit}, after taking $\liminf_{x\rightarrow 0}$ we obtain
\[
    \liminf_{x\rightarrow0}\frac{\vert x\vert^{r-1}\ln (n(k+1))!}{\ln\vert x\vert} \geq  -\frac{k+1}{\alpha}.
\]

For the upper limit,
\[
    \ln(n(k+1))! \geq \int_{1}^n\ln s\; ds= n(\ln n-1)+1.
\]
Using analogous arguments, one gets 
\[
    \limsup_{x\rightarrow0}\frac{\vert x\vert^{r-1}\ln (n(k+1))!}{\ln\vert x\vert} \leq  -\frac{k+1}{\alpha},
\]
and the result follows.

\end{proof}

We now proceed to prove Theorem~\ref{THMasymptoticNonhyperbolic}. Let $\varepsilon\in(0,\alpha)$ be sufficiently small and $x_0<0$ close to $0$, so that conditions C1--C5 hold.

\begin{proof}[Proof of Theorem \ref{THMasymptoticNonhyperbolic}]
	We again divide the proof into two parts, for the upper and lower limits, respectively.
    \medskip
	
	\noindent\textit{Part 1: lower bound in (\ref{ordernonhyperbolic}).}
    
    \noindent For $n=n_{x_0}^x$, we employ the upper bound obtained in \eqref{ineqMainLemma_Nonhyperbolic} from the nonhyperbolic case of Lemma~\ref{LEMMAcrucial}. By taking logarithm, it yields
\begin{align*}
    \ln\phi(x) &\leq n \ln(C_2\cdot k!)-
    \ln\; (n(k+1))! + n(k+1)\ln\vert X_{-1}\vert\\
    & \qquad -(k+1)\sum_{m=1}^{n-1}\sum_{j=1}^m
    \ln\left(h_+'(X_{-j-1})\right) + \ln B.
\end{align*}
The first and the last term in the right hand side will vanish after multiplying times $\vert x\vert^{r-1}/\ln\vert x\vert$ and taking limits as $x\rightarrow0$. Therefore, 
\begin{equation}
\label{eq:last_nonhyperbolic}
    \liminf_{x\rightarrow 0}\frac{\vert x \vert^{r-1} \ln\phi(x)}{\ln\vert x\vert} \geq \frac{r(k+1)}{\alpha(r-1)} - (k+1)\lim_{x\rightarrow 0} \sum_{m=1}^{n-1}\sum_{j=1}^m
    \ln\left(h_+'(X_{-j-1})\right),
\end{equation}
if the limit on the right hand side exists. This is indeed the case, as we show in the following.

\smallskip
\noindent\textit{Claim: the limit on the right hand of \eqref{eq:last_nonhyperbolic} side vanishes.} Indeed, from \eqref{MapNonHyperbolic} it follows that for all $x$ sufficiently close to $0$,
\begin{equation}
\label{eq:bounds_derivative}
    1-\alpha_1^\varepsilon r\vert x\vert^{r-1}\leq h_+'(x) \leq 1-\alpha_2^\varepsilon r\vert x\vert^{r-1},
\end{equation}
 where $\alpha^\varepsilon_1:=\alpha+\varepsilon$ and $\alpha^\varepsilon_2:=\alpha- \varepsilon$. Denote the double sum as $S$, so that
 \begin{align*}
     S & :=\sum_{m=1}^{n-1}\sum_{j=1}^m\ln \left( h_+'(X_{-j-1}) \right) \geq
     \sum_{m=1}^{n}\sum_{j=1}^m\ln \left( h_+'(X_{-j})\right) \geq 
     \sum_{m=1}^{n}\sum_{j=1}^m\ln \left( 1-\alpha_1^{\varepsilon} r\vert X_{-j}\vert^{r-1}\right) \\
     &=
     \sum_{m=1}^{n}\sum_{j=1}^m\frac{\ln \left( 1-\alpha_1^{\varepsilon} r\vert X_{-j}\vert^{r-1} \right)}{X_{-j+1}-X_{-j}}\Delta X_{-j},
\end{align*}
where $\Delta X_{-j}:= X_{-j+1}-X_{-j}$.

From \eqref{MapNonHyperbolic}, by taking $x_0$ sufficiently close to $0$, we have that for all $x\in[h_+^{-1}(x_0),0]$, 
\[
    x+\alpha_2^\varepsilon \vert x\vert^r\leq h_+(x)\leq x+\alpha_1^\varepsilon \vert x\vert^r.
\]
In particular, $X_{-j+1}-X_{-j}=h_+\left(  X_{-j}\right)-X_{-j}\geq \alpha_2^\varepsilon \vert X_{-j}\vert^r$, 
so that
\[
   \frac{\ln(1-\alpha_1^\varepsilon\vert X_{-j}\vert^{r-1})}{X_{-j+1}-X_{-j}} \geq \frac{\ln(1-\alpha_1^\varepsilon\vert X_{-j}\vert^{r-1})}{\alpha_2^\varepsilon\vert X_{-j}\vert^r}.
\]
Hence,
\begin{equation}
\label{eq:lower_1}
    S\geq \sum_{m=1}^{n}\sum_{j=1}^m\frac{\ln \left( 1-\alpha_1^{\varepsilon} r\vert X_{-j}\vert^{r-1}\right)}{\alpha_2^{\varepsilon}\vert X_{-j}\vert^r}\Delta X_{-j}.
\end{equation}

Notice that the lower bound in \eqref{eq:lower_1} resembles a Riemann sum. We further estimate $S$ from below with an integral. Consider the function 
\[
f(s)=\frac{\ln(1-\alpha_1^\varepsilon rs^{r-1})}{\alpha_2^\varepsilon s^{r}}, \qquad s>0.
\]
Notice that $f(s)\rightarrow -\infty$ as $s\rightarrow0$ and that for $s$ sufficiently small, the function $f$ is increasing. The inner sum is thus bounded from below as
\[    \sum_{j=1}^m \frac{\ln(1-\alpha_1^{\varepsilon}r \vert X_{-j}\vert^{r-1})}{\alpha_2^\varepsilon \vert X_{-j}\vert^r}\Delta X_{-j} 
    \geq \int_{\vert x\vert}^{\vert X_{-m}\vert} \frac{\ln(1-\alpha_1^\varepsilon r s^{r-1})}{\alpha_2^\varepsilon s^r}ds.
\]
We can further estimate $S$ from below, continuing from \eqref{eq:lower_1}, as
\begin{equation}
\label{eq:lower_2}
 S\geq \sum_{m=1}^{n} \int_{\vert x\vert}^{\vert X_{-m}\vert}\frac{\ln(1-\alpha_1^{\varepsilon} r s^{r-1})}{\alpha_2^{\varepsilon} s^r}ds.
\end{equation}
We refer the reader to Figure~\ref{fig:integrals_approx} (a) where the argument above is graphically summarised.

Recall that $\ln(1-x)\geq -2x$ for positive $x\approx 0$. By using this inequality in \eqref{eq:lower_2}, we obtain
\begin{align*}
     S & \geq -\frac{2r\alpha_1^\varepsilon}{\alpha_2^\varepsilon} \sum_{m=1}^{n}\int_{\vert x\vert}^{\vert X_{-m}\vert}\frac{ds}{s} = -\frac{2r\alpha_1^\varepsilon}{\alpha_2^\varepsilon}\sum_{m=1}^{n}\left( 
    \ln\vert X_{-m}\vert - \ln\vert x\vert\right)  \\
    & = -\frac{2r\alpha_1^\varepsilon}{\alpha_2^\varepsilon}\sum_{m=1}^{n} \ln\vert X_{-m}\vert+ \frac{2r\alpha_1^\varepsilon}{\alpha_2^\varepsilon}(n-1)\ln\vert x\vert\\
    & = -\frac{2r\alpha_1^\varepsilon}{\alpha_2^\varepsilon}\sum_{m=1}^{n} \frac{\ln\vert X_{-m}\vert}{X_{-m}-X_{-m-1}}\Delta X_{-m-1} + \frac{2r\alpha_1^\varepsilon}{\alpha_2^\varepsilon}(n-1)\ln\vert x\vert,
\end{align*}
where $\Delta X_{-m-1}=X_{-m}-X_{-m-1}$. Since $h_+^{-1}(x)=x-\alpha \vert x\vert^r+o(\vert x\vert^r)$, if $x_0$ is sufficiently close to $0$, we have that $X_{-m}-X_{-m-1}\leq \alpha_1^\varepsilon \vert X_{-m}\vert^r$ for all $m=1,\ldots,n$. It follows that
\begin{equation}
\label{eq:lower_3}
    S\geq -\frac{2r}{\alpha_2^{\varepsilon}}\sum_{m=1}^{n}
    \frac{\ln\vert X_{-m}\vert}{ \vert X_{-m}\vert^r}\Delta X_{-m-1}
    + \frac{2r\alpha_1^\varepsilon}{\alpha_2^\varepsilon}(n-1)\ln\vert x\vert.
\end{equation}
Similarly as before, consider the function 
\[
    g(s)=\frac{\ln s}{s^r}, \qquad s>0.
\]
Notice that $g$ is increasing for $s\approx0$ and $g(s)\rightarrow-\infty$ as $s\rightarrow0$. Therefore, from \eqref{eq:lower_3}, we estimate $S$ from below with an integral as (see Figure~\ref{fig:integrals_approx} (b))
\begin{equation}
\label{eq:lower_4}
    S\geq -\frac{2r}{\alpha_2^{\varepsilon}} \int_{\vert X_{-1}\vert}^{\vert X_{-n-1}\vert}\frac{\ln s}{s^r}ds + \frac{2r\alpha_1^\varepsilon}{\alpha_2^\varepsilon}(n-1)\ln\vert x\vert.
\end{equation}

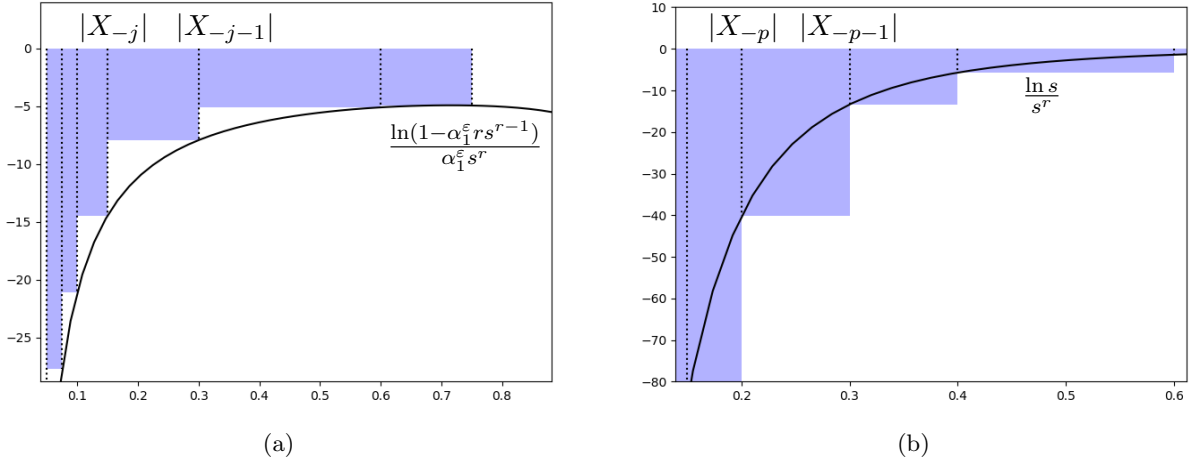
\begin{figure}
     \centering
     \begin{subfigure}[b]{0.47\textwidth}
         \centering
          \begin{overpic}[width=\linewidth]{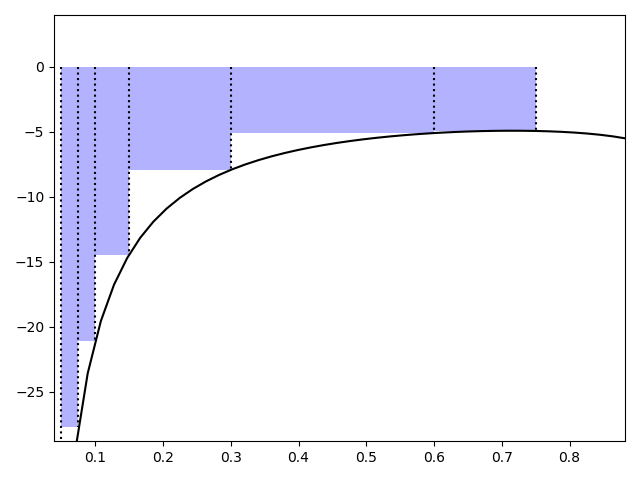}
          \put(60,46){$f(s)=\frac{\ln(1-\alpha_1^\varepsilon r s^{r-1})}{\alpha_1^{\varepsilon}s^r}$}
          \put(32,67){$\vert X_{-j}\vert$}
          \put(15,67){$\vert X_{-j+1}\vert$}
      \end{overpic}
      \caption{}
         \label{subfig:integral_comp1}
     \end{subfigure}
     \hfill
     \begin{subfigure}[b]{0.47\textwidth}
         \centering
          \begin{overpic}[width=\linewidth]{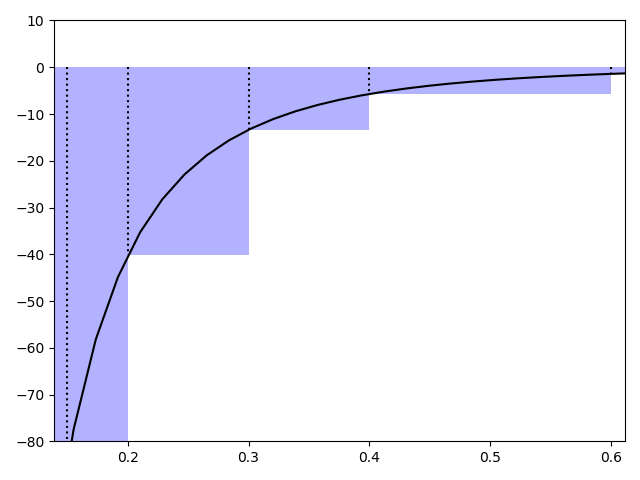}
           \put(69,55){$g(s)=\frac{\ln s}{s^r}$}
           \put(14,67){$\vert X_{-m}\vert$}
           \put(30,67){$\vert X_{-m-1}\vert$}
      \end{overpic}
         \caption{}
    \label{subfig:integral_comp2}
     \end{subfigure}
     \hfill
       \caption[Integral comparisons]{In (a) and (b), a graphical summary of the integral estimation of the sums \eqref{eq:lower_1} and \eqref{eq:lower_3} are respectively sketched.}
        \label{fig:integrals_approx}
\end{figure}

By solving the indefinite integral using integration by parts (taking $u=\ln s$ and $dv=s^{-r}ds$) we have that
\[
    \int\frac{\ln s}{s^r}ds =
    \frac{1}{(r-1)s^{r-1}}\left[ \frac{1}{r-1} -\ln s \right].
\]
Since $\vert X_{-n-1}\vert \in [h_+^{-2}(x_0),h_+^{-1}(x_0)]$, there is a constant $E$ such that
\begin{equation}
    \label{eq:lower_5}
    S\geq \frac{2r}{\alpha_2^{\varepsilon}(r-1)\vert X_{-1}\vert^{r-1}}\left( \frac{1}{r-1} -\ln \vert X_{-1}\vert \right) + \frac{2r\alpha_1^\varepsilon}{\alpha_2^\varepsilon}(n-1)\ln\vert x\vert + E.
\end{equation}
Observe that
\[
    \frac{\vert x\vert}{\vert X_{-1}\vert}\rightarrow 1, \qquad \frac{\ln \vert x\vert}{\ln \vert X_{-1}\vert}\rightarrow 1, \qquad \textup{as } x\rightarrow 0.
\]
Hence, after multiplying times $\vert x\vert^{r-1}/\ln\vert x\vert<0$ on both sides, and taking $\limsup$ as $x\rightarrow0$, using Lemma~\ref{LEMMAscalingdeterministic} one gets
\[
    0<\limsup_{x\rightarrow0} \frac{\vert x\vert^{r-1} S}{\ln\vert x\vert}\leq \frac{2r}{\alpha(r-1)}-\frac{2r}{\alpha_2^\varepsilon(r-1)}.
\]
Since it holds for all $\varepsilon>0$, and $\alpha_2^{\varepsilon}=\alpha- \varepsilon$, the claim follows.

Finally, since the limit in the right hand side of \eqref{eq:last_nonhyperbolic} vanishes, we can conclude that
\[
    \liminf_{x\rightarrow0} \frac{\vert x\vert \ln\phi(x)}{\ln\vert x\vert}\geq 
    \frac{r(k+1)}{\alpha(r-1)}.
\]

\noindent\textit{Part 2: upper bound in (\ref{ordernonhyperbolic}).} 

\noindent We first provide a rough lower bound for $\phi(x)$, by reducing the integration interval in the lower bound in \eqref{ineqMainLemma_Nonhyperbolic} to $[X_{-1},x]$, yielding
\begin{align*}
    \phi(x) & \geq \frac{(C_1 k!)^n}{[n(k+1)-1]!}\int_{X_{-1}}^x (y-X_{-1})^{n(k+1)-1} \phi(Y_{-n+1})dy \\
    & \geq \frac{(C_1 k!)^n}{[n(k+1)]!}(x-X_{-1})^{n(k+1)} \cdot D \geq 
    \frac{\tilde{C}_1^n}{[n(k+1)]!}\vert x \vert^{rn(k+1)} \cdot D,
\end{align*}
where $\tilde{C}_1:=(\alpha-\varepsilon)^{k+1}C_1 k!$ and $D=\inf_{x\in[x_0,h_+(x_0)]}\phi(x)$. By taking $\ln$ in the last expression, and since $X_{-n+1}\in(x_0,h_+(x_0))$, we obtain
\begin{align*}
    \ln\phi(x)& \geq  rn(k+1)\ln\vert x\vert-\ln(n(k+1)!)
+n\ln\tilde{C}_1+\ln D.
\end{align*}
By multiplying times $\frac{\vert x\vert^{r-1}}{\ln\vert x\vert}<0$, and using Proposition~\ref{LEMMAscalingdeterministic} and Corollary~\ref{CORO: log_factorial_nonhyp}, it follows that
\begin{equation}
\label{eq:rough}
    \limsup_{x\rightarrow 0}\frac{\vert x\vert^{r-1} \ln\phi(x)}{\ln\vert x\vert}\leq \frac{r(k+1)}{\alpha(r-1)} +\frac{(k+1)}{\alpha} = \frac{(2r-1)(k+1)}{\alpha(r-1)}.
\end{equation}
    
While this is still not the value we aim to obtain, we use this result to state that for any $\delta>0$, there exists $\tilde{x}_0<0$ such that for all $x\in(\tilde{x}_0,0)$
\begin{equation}
    \label{eq:estimate_log_phi}
    \ln\phi(x)\geq \delta \left(\frac{\ln \vert x\vert}{\vert x\vert^r}\right).
\end{equation}
This is indeed the case since 
\[
    \frac{\vert x\vert^{r}\ln \phi(x)}{\ln\vert x\vert}\rightarrow 0, \qquad \textup{as } x\rightarrow0
\]

We are now ready to improve \eqref{eq:rough} and conclude the proof. By taking $\ln$ directly in \eqref{ineqMainLemma_Nonhyperbolic}, by means of Jensen's inequality, we get
\begin{equation}
\label{eq:separation_final}
\begin{split}
    \ln\phi(x) & \geq n\ln(C_1k!)-\ln(n(k+1)-1)!+\ln\vert X_{-1}\vert   \\
    & \qquad +\frac{n(k+1)}{\vert X_{-1}\vert}\int_{X_{-1}}^0\ln(y-X_{-1}) dy + \frac{1}{\vert X_{-1}\vert}\int_{X_{-1}}^0 \ln\phi(Y_{-n+1})dy,
\end{split}
\end{equation}
The first integral is explicitly solvable, so that
\begin{equation}
    \label{eq:first_integral}
\frac{n(k+1)}{\vert X_{-1}\vert}\int_{X_{-1}}^0\ln(y-X_{-1}) dy =
n(k+1)\left( \ln\vert X_{-1}\vert -1 \right).
\end{equation}
For the second integral, we use \eqref{eq:estimate_log_phi} to show the following statement. 

\noindent\textit{Claim: there exists $M\equiv M(x_0)$ such that 
\begin{equation}
\label{eq:second_integral}
    \int_{X_{-1}}^0\ln\phi(Y_{-n+1})dy \geq M.
\end{equation}}
Indeed, assume without loss of generality that $x_0$ is sufficiently close to $0$ so that \eqref{eq:estimate_log_phi} holds for $\delta=1$ and all $x\in(h_+^{-1}(x_0),0)$. Therefore, by doing the change of variable $z=Y_{-n+1}$,
\begin{align*}
    \int_{X_{-1}}^0\ln\phi(Y_{-n+1}) & = \int_{X_{-n}}^0{h'(z)\cdots h'(Z_{n-1})\ln\phi(z)dz} \geq \int_{X_{-n}}^0\ln\phi(z)\geq \int_{h_+^{-1}(x_0)}^0\ln\phi(z)dz \\
    &\geq \int_0^{\vert h_+^{-1}(x_0)\vert} \frac{\ln s}{s^r} = \frac{\vert h_+^{-1}(x_0)\vert}{r-1}\left[ (r-1)\ln\vert h_+^{-1}(x_0)\vert -1\right] \equiv M(x_0).
\end{align*}

Combining \eqref{eq:first_integral} and \eqref{eq:second_integral} in \eqref{eq:separation_final}, we get
\[
    \ln\phi(x)\geq n\ln(C_1 k!)-\ln (n(k+1)-1)! + \ln\vert X_{-1}\vert + n(k+1)(\ln\vert X_{-1}-1\vert) + \frac{M}{\vert X_{-1}\vert}.
\]
We conclude that
\[
    \frac{\vert x\vert^{r-1}\ln\phi(x)}{\ln\vert x\vert} \leq - \frac{\vert x\vert^{r-1}\ln (n(k+1))!}{\ln \vert x\vert}+ (\vert x\vert^{r-1} \cdot n) (k+1) + o(1).
\]

It follows from Lemma~\ref{LEMMAscalingdeterministic} and Corollary~\ref{CORO: log_factorial_nonhyp} that
\[
    \limsup_{x\rightarrow 0} \frac{\vert x\vert^{r-1}\ln\phi(x)}{\ln\vert x\vert}\leq \frac{k+1}{\alpha} + \frac{k+1}{\alpha(r-1)}=\frac{(k+1)r}{\alpha(r-1)},
\]
concluding Step 2. Combined with Step 1, the theorem follows.
\end{proof}

Let $d\mu=\phi\; d\Leb$ be the unique stationary distribution for (\ref{RDE}) supported on $M$. Consider the tail of $\mu$, that is the function $\mathcal{T}(x):=\mu\left( [x,x_+] \right)$. We now provide an analogous result to Theorem~\ref{THMasymptoticNonhyperbolic} for the tail $\mathcal{T}$.

\begin{coro}
	\label{CORO:tail_Measure_nonhyper}
		Assume that the hypothesis of Theorem~\ref{THMasymptoticNonhyperbolic} hold. Then,
		\begin{equation}
		\label{LDPtailsNH}
		\lim_{x\rightarrow x_+}\frac{(x_+-x)^{r-1}\ln\mathcal{T}(x)}{\ln(x_+-x)}= \frac{r(k+1)}{\alpha(r-1)}.
		\end{equation}
\end{coro}
\begin{proof}
	The proof follows analogous arguments as those for Corollary~\ref{CORO:tail_measure_hyperbolic}.
\end{proof}

\bigskip

\subsection*{Acknowledgments}
 The authors would like to thank Dmitry Turaev, Sebastian van Strien, Ale Jan Homburg, Andrew Clarke, and Maximilian Engel for fruitful and insightful discussions during the early stage of this research. They also thank Wei Hao Tey for pointing out some inaccuracies in an earlier version of the manuscript, and for discussions on the practical applicability of the results, which derived into the work \cite{OliconTey25}. 

\subsection*{Funding}
 G.O.-M. thanks CONAHCYT (formerly CONACYT) and the Department of Mathematics of Imperial College London for a PhD scholarship. This paper builds on Chapter~3 of the PhD thesis~\cite{Olicon21}. GOM also thanks the Deutsche Forschungsgemeinschaft (DFG) through grant CRC 1114 ``Scaling Cascades in Complex Systems'', Project Number 235221301, Project A02 ``Multiscale data and asymptotic model assimilation for atmospheric flows '', and through Germany's Excellence Strategy -- The Berlin Mathematics Research Center MATH+ EXC-2046/1, project 390685689,  subproject AA1-8.

 J.S.W.L. and M.R. are grateful for support by the EPSRC research grants EP/W009455/1 and EP/Y020669/1. JSWL gratefully acknowledges support from the EPSRC Centre for Doctoral Training in Mathematics of Random Systems: Analysis, Modelling and Simulation (EP/S023925/1), and thanks IRCN (Tokyo) and GUST (Kuwait) for their research support.


\bibliographystyle{plain} 

\bibliography{biblio}       


\end{document}